\documentclass[12pt]{article}
\usepackage{amsfonts}
\usepackage{amssymb}
\usepackage{amsmath}
\usepackage{graphicx}
\usepackage{subfigure}
\usepackage{theorem}

\def\date{3 April 2013}  
\newtheorem{proposition}{proposition}[section]
\newtheorem{lemma}[proposition]{Lemma}
\newtheorem{corollary}[proposition]{Corollary}
\newtheorem{theorem}[proposition]{Theorem}

\newtheorem{conjecture}[proposition]{Conjecture}
\newcommand{\qed}{$\square$\bigskip}

\newcommand{\proof}{{\noindent\bf Proof. }}
{\theorembodyfont{\rmfamily}

\newtheorem{hypothesis}[proposition]{Hypothesis}
}
\def\mytextindent#1{\indent\indent\llap{\rm#1\enspace}\ignorespaces}
\def\myitem{\par\hangindent0pt\mytextindent}

\def\myclaim#1#2\par{{\medbreak\noindent\rlap{\rm(#1)}\ignorespaces
 \rightskip20pt
 \hangindent=20pt\hskip20pt{\ignorespaces\sl#2}\smallskip}}
\def\junk#1{}

\def\mylabel#1{{\label{#1}\marginpar{#1}}}

\font\smallrm=cmr8
\font\smallrm=cmr8
\begin{document}
\newif\ifproofmode
\baselineskip=12pt
\phantom{a}\vskip .25in
\centerline{{\bf  $K_6$ MINORS IN 6-CONNECTED GRAPHS OF BOUNDED TREE-WIDTH}}
\vskip.4in
\centerline{\bf Ken-ichi Kawarabayashi\footnote{\smallrm Supported by
JST, ERATO, Kawarabayashi Large Graph Project}}
\centerline{National Institute of Informatics}
\centerline{2-1-2 Hitotsubashi, Chiyoda-ku, Tokyo 101-8430, Japan}
\medskip
\centerline{{\bf Serguei Norine}\footnote{\smallrm
Partially supported by NSF under Grants No.~DMS-0200595 and~DMS-0701033. }}
\centerline{Department of Mathematics}
\centerline{Princeton University}
\centerline{Princeton, NJ 08544, USA}
\medskip
\centerline{{\bf Robin Thomas}%
\footnote{\smallrm Partially supported by NSF under
Grants No.~DMS-0200595, DMS-0354742, and~DMS-0701077.}}
\centerline{School of Mathematics}
\centerline{Georgia Institute of Technology}
\centerline{Atlanta, Georgia  30332-0160, USA}
\medskip
\centerline{and}
\medskip
\centerline{\bf Paul Wollan\footnote{\smallrm Partially supported by the European Research Council under the European Union’s Seventh Framework Programme (FP7/2007-2013)/ERC Grant Agreement no. 279558.}}
\centerline{Department of Computer Science}
\centerline{University of Rome ``La Sapienza"}
\centerline{00185 Rome, Italy}

\vskip 1in \centerline{\bf ABSTRACT}
\bigskip
\parshape=1.5truein 5.5truein
We prove that every sufficiently large $6$-connected graph of
bounded tree-width either has a $K_6$ minor, or has
a vertex whose deletion makes the graph planar.
This is a step toward proving that the same conclusion holds
for all sufficiently large $6$-connected graphs.
J\o rgensen conjectured that it holds for all $6$-connected graphs.

\vfill \baselineskip 11pt \noindent 8 April 2005, revised \date.
\vfil\eject

\ifproofmode
  \def\sectionbreak{\vfil\eject}
   \newcount\remarkno
   \def\REMARK#1{{%
      \footnote{\baselineskip=11pt #1
      \vskip-\baselineskip}\global\advance\remarkno by1}}
    \def\mylabel#1{{\label{#1}\marginpar{#1}}}

\else
  \def\sectionbreak{}
  \def\REMARK#1{}
  \def\mylabel#1{\label{#1}}
\fi

\baselineskip 18pt

\section{Introduction}\mylabel{sec:intro}

Graphs in this paper are allowed to have loops and multiple edges.
A graph is a {\em minor} of another if the first can be obtained from
a subgraph of the second by contracting edges. An {\em $H$ minor}
is a minor isomorphic to $H$.
A graph $G$ is {\em apex} if it has a vertex $v$ such that
$G\backslash v$ is planar.
(We use $\backslash$ for deletion.) J\o rgensen~\cite{Jor} made the following beautiful conjecture.

\begin{conjecture}
\mylabel{con:jorgensen}
Every $6$-connected graph with no $K_6$ minor is apex.
\end{conjecture}

In a companion paper~\cite{KawNorThoWollarge} 
we prove that Conjecture~\ref{con:jorgensen}
holds for all sufficiently big $6$-connected graphs. Here we
establish the first step toward that goal, the following.

\begin{theorem}
\mylabel{main}
For every integer $w\ge1$ there exists an integer $N$ such that
every $6$-connected graph on at least $N$ vertices and tree-width
at most $w$ with no $K_6$ minor is apex.
\end{theorem}

We define tree-width later in this section, but let us discuss
J\o rgensen's conjecture first.
It is related to Hadwiger's conjecture~\cite{Had}, the following.

\begin{conjecture}
\mylabel{con:hadwiger}
For every integer $t\ge1$, if a loopless graph has no $K_t$ minor, then it
is $(t-1)$-colorable.
\end{conjecture}

Hadwiger's conjecture is known for $t\le6$.
It is trivial for $t\le3$, and is still fairly easy for $t=4$,
as shown by Dirac~\cite{Dirproperty}.
However, for $t\ge5$ Hadwiger's conjecture implies the Four-Color Theorem.
Wagner~\cite{Wag37} gave a structural characterization of graphs
with no $K_5$ minor, which implies that  Hadwiger's conjecture for $t=5$
is actually equivalent to the Four-Color Theorem.
The same conclusion has been obtained for $t=6$ in~\cite{RobSeyThoHad} by showing
that a minimal counterexample to Hadwiger's conjecture for $t=6$
is apex.
The proof uses an earlier result of Mader~\cite{MadHomkrit}
that every minimal counterexample to Conjecture~\ref{con:hadwiger}
is $6$-connected.
Thus Conjecture~\ref{con:jorgensen}, if true, would give more
structural information.
Furthermore, the structure of all graphs with no $K_6$ minor is not known,
and appears complicated and difficult.
Thus obtaining a structural characterization of graphs with no
$K_6$ minor, an analogue of Wagner's theorem mentioned above,
 appears beyond reach at the moment.
On the other hand, Conjecture~\ref{con:jorgensen} provides a nice
necessary and sufficient condition for $6$-connected graphs.
Unfortunately, it, too, appears to be a difficult problem.

Let us turn to tree-width and our proof method.
Tree-width of a graph was first defined by Halin~\cite{HalSfunctions},
and was later rediscovered in~\cite{RobSeyGM3}, and, independently,
in~\cite{ArnProLinear}. The definition is as follows.
A {\em tree-decomposition} of a graph $G$ is a pair $(T,Y)$, where $T$
is a tree and $Y$ is a family $\{Y_t \mid t \in V(T)\}$ of vertex sets
$Y_t\subseteq V(G)$, such that the following two properties hold:

\begin{enumerate}
\item[(W1)] $\bigcup_{t \in V(T)} Y_t = V(G)$, and every edge of $G$ has
both ends in some $Y_t$.
\item[(W2)] If  $t,t',t''\in V(T)$ and $t'$ lies on the path in $T$
between $t$ and $t''$, then $Y_t \cap Y_{t''} \subseteq Y_{t'}$.
\end{enumerate}
The {\em width} of a tree-decomposition $(T,Y)$ is $\max_{t\in V(T)} (|Y_t|-1)$,
and the {\em tree-width} of a graph $G$ is the minimum width of a
tree-decomposition of $G$.

Our proof of Theorem~\ref{main} proceeds as follows. We choose a tree-decomposition $(T,W)$ of
$G$ of width $w$ with no ``redundancies".
It follows easily that if $T$ has a vertex of large degree, then
$G$ has a $K_6$ minor, and so we may assume that $T$ has a long path.
For the rest of the proof we concentrate our effort on this long path.
Since other branches of $T$ are inconsequential, we convert $(T,W)$
to a ``linear decomposition", which is really just a tree-decomposition,
where the underlying tree is a path, but we find it more convenient
to number the sets of vertices $W_0,W_1,\ldots,W_l$, rather than
index them by the vertices of a path. At this point we no longer
require that the width be bounded; all we need is that the intersections
$W_{i-1}\cap W_i$ are bounded and that $l$ is sufficiently large.
Thus we may assume (by trimming our linear decomposition) that
all the sets $W_{i-1}\cap W_i$ have the same size, say $q$. Furthermore,
it can be arranged (by invoking the result from~\cite{RThoMenger}
or by a direct argument) that there exist $q$ disjoint paths
$P_1,P_2,\ldots,P_q$
from $W_{0}\cap W_1$ to $W_{l-1}\cap W_l$.
We apply the pigeon hole principle many times, each time trimming
the linear decomposition, but still keeping it sufficiently long,
to make sure that if the subgraph $G[W_i]$ has some useful property for some
$i\in\{1,2,\ldots,l-1\}$, then all the graphs $G[W_i]$ have that property
for all $i\in\{1,2,\ldots,l-1\}$.

A prime example of a useful property is the existence of
two disjoint paths $Q_1,Q_2$ in $G[W_i]$,
internally disjoint from $P_1,P_2,\ldots,P_q$,
 with ends $u_1,v_1$ and $u_2,v_2$, respectively, such that
$u_1,v_2\in V(P_1)$, $u_2,v_1\in V(P_2)$ and they appear on those
paths in the order listed as $P_1$ and $P_2$ are traversed from
$W_{0}\cap W_1$ to $W_{l-1}\cap W_l$.
In those circumstances we say that $P_1$ and $P_2$ {\em twist} in $W_i$.
Thus, in particular, we can arrange that if two paths $P_j$ and $P_{j'}$ twist in
$W_i$ for some $i\in\{1,2,\ldots,l-1\}$, then they twist in
$W_i$ for all $i\in\{1,2,\ldots,l-1\}$.
On the other hand, if two paths $P_j$ and $P_{j'}$  twist in
$W_i$ for all $i\in\{1,2,\ldots,l-1\}$ and $l$ is not too small, then
$G$ has a $K_6$ minor.
This is the sort of argument we will be using, but the details are too
numerous to be described in their entirety here.

In~\cite{KawNorThoWollarge} we use Theorem~\ref{main} to prove
J\o rgensen's conjecture for sufficiently big graphs, formally the
following:

\begin{theorem}
\mylabel{future}
There exists an integer $N$ such that
every $6$-connected graph on at least $N$ vertices  with no $K_6$ minor is apex.
\end{theorem}

How does Theorem~\ref{main} help in the proof of Theorem~\ref{future}?
By the excluded grid theorem of Robertson and Seymour~\cite{RobSeyGM5}
(see also \cite{DieGorJenTho,ReedBCC,RobSeyThoQuickPlanar})
\nocite{DieGorJenTho}
\nocite{ReedBCC}
\nocite{RobSeyThoQuickPlanar}
it suffices to prove Theorem~\ref{future} for graphs that have a sufficiently large
grid minor. We then analyze how the remainder of the graph attaches to
the grid. We refer to~\cite{KawNorThoWollarge} for details.

The paper is organized as follows.
In Section~\ref{sec:reroute} we state a few lemmas, mostly from other papers.
In Section~\ref{sec:linear} we convert a tree-decomposition into a
linear decomposition,
as described above, and we prove that the linear decomposition can be
chosen with some additional useful properties.
In Section~\ref{sec:auxil} we introduce the auxiliary graph---its vertices
are the paths $P_1,P_2,\ldots,P_q$, and two of them are adjacent if
they are joined by a path avoiding all the other paths $P_1,P_2,\ldots,P_q$.
By joined we mean in some or every $W_i$; by now the two are equivalent.
We use the auxiliary graph to further refine the linear decomposition.
A {\em core} is a component of the subgraph of the auxiliary graph induced
by those of the paths  $P_1,P_2,\ldots,P_q$ that have at least one edge.
We show, among other things, that every core is a path or a cycle.
In Section~\ref{sec:pin} we use the theory of ``non-planar extensions"
of planar graphs from~\cite{RobSeyThoExt} to get under control adjacencies
in the auxiliary graph of those paths $P_i$ that are trivial.
In Section~\ref{sec:taming} we further refine our linear decomposition
to arrange that the part of $G$ that corresponds to a core can be drawn
either in a disk or in a cylinder, depending on whether the core is
a path or a cycle.
In Section~\ref{sec:control} we digress and prove a slight extension
of a result of DeVos and Seymour~\cite{DevSeyExt3col}.
Finally, in Section~\ref{sec:cyl} we essentially complete the proof of Theorem~\ref{main} in the case when some
core of the auxiliary graph is a cycle, and in Section~\ref{sec:planar}
we do the same when some core is a path.


\section{Rerouting and rural societies}\mylabel{sec:reroute}

Let $S$ be a subgraph of a graph $G$.  An {\em $S$-bridge} in $G$ is a
connected subgraph $B$ of $G$ such that $E(B)\cap E(S)=\emptyset$
and either $E(B)$ consists of a unique edge with both ends in $S$, or for
some component $C$ of $G\backslash V(S)$ the set $E(B)$ consists of all
edges of $G$ with at least one end in $V(C)$.  The vertices in
$V(B)\cap V(S)$ are called the {\em attachments} of $B$.  We say that an $S$-bridge $B$ \emph{attaches to} a subgraph $H$ of $S$ if $V(H) \cap V(B) \neq \emptyset$.

Now let $S$ be such that no block of $S$ is  a cycle.
By a {\em segment of $S$} we mean a maximal subpath $P$ of $S$ such that
every internal vertex of $P$ has degree two in $S$.
It follows that the segments 
of $S$ are uniquely determined.  Now if $B$ is an $S$-bridge of $G$, then we
say that $B$ is {\em unstable} if some segment of $S$ includes all the
 attachments of $B$, and otherwise we say that $B$ is
{\em stable}.  Our next lemma says that it is possible to make all
bridges stable by making the following ``local" changes.  Let $G$ and $S$ be
as before, let $P$ be a segment of $S$ of length at least two,
and let $Q$ be a path
in $G$ with ends $x,y\in V(P)$ and otherwise disjoint from $S$.
Let $S'$ be obtained from $S$ by replacing the path $xPy$
(the subpath of $P$ with ends $x$ and $y$) by $Q$;
then we say that
$S'$ was obtained from $S$ by {\em rerouting} $P$ along $Q$, or
simply that $S'$ was obtained from $S$ by {\em rerouting}.  Please
note that $P$ is required to have length at least two,
and hence this relation is not symmetric.  We say that the
rerouting is {\em proper} if
all the attachments
of the $S$-bridge that contains $Q$  belong to $P$.
The following lemma is essentially due to Tutte.

\begin{lemma}
\mylabel{prestable}
Let $G$ be a simple graph,
and let $S$ be a subgraph of $G$
such that no block of $S$ is a cycle.
Then there exists a subgraph $S'$ of $G$
 obtained from $S$ by a sequence of proper reroutings
such that if an $S'$-bridge $B$ of $G$ is unstable, say all its attachments
belong to a segment $P$ of $S'$, then there exist
vertices
$x,y\in V(P)$ such that
some component of $G\backslash\{x,y\}$ includes a vertex of $B$
and is disjoint from $S'\backslash V(P)$.
\end{lemma}

\proof  We may choose a subgraph $S'$ of $G$ obtained from $S$ by a
sequence of proper reroutings such that the number of vertices that
belong to stable $S'$-bridges is maximum, and, subject to that,
$|V(S')|$ is minimum.  We will show that $S'$ is as
desired.  To that end we may assume that
$B$ is an $S'$-bridge of $G$ with all its attachments in a
segment $P$ of $S'$.

Let $v_0, v_1,\dots, v_k$ be distinct vertices of $P$, listed in order of
occurrence on $P$ such that $v_0$ and $v_k$ are the ends of $P$ and
$\{v_1,\dots, v_{k-1}\}$ is the set of all internal vertices of $P$
that are attachments of a stable $S'$-bridge.
We claim that if
$u,v$ are two attachments of $B$, then no $v_i$ belongs to
the interior of $uPv$.
To prove this suppose to the contrary that $v_i$ is an internal
vertex of $uPv$.
But then replacing $uPv$ by an induced subpath of $B$ with ends $u,v$
and otherwise disjoint from $S'$ is a proper rerouting that produces
a graph $S''$ with strictly more vertices belonging to stable
$S''$-bridges, contrary to the choice of $S'$.
This proves our claim that no $v_i$ belongs to the interior of $uPv$.
But then for some $i=1,2,\ldots,k$ the path $v_{i-1}Pv_i$
includes all attachments of $B$.
Since $G$ has no parallel edges, the same argument (using the minimality
of $|V(S')|$) now shows that $V(B)-\{v_{i-1},v_i\} \neq \emptyset$.
Consequently some component $J$ of $G\backslash \{v_{i-1},v_i\}$
includes a vertex of $B$. It follows that $B\backslash \{v_{i-1},v_i\}$
is a subgraph of $J$.
Now $B$ has all its attachments in
$v_{i-1}Pv_i$, the interior of $v_{i-1}Pv_i$ includes no
attachment of a stable $S'$-bridge, and
(by what we have shown about $B$) every unstable $S'$-bridge with
an attachment in the interior of $v_{i-1}Pv_i$ has all its attachments
in  $v_{i-1}Pv_i$.
It follows that $J$ is disjoint
from $S'\backslash V(P)$, as desired.~\qed


We deduce the following corollary.

\begin{theorem}
\mylabel{stable}
Let $G$ be a $3$-connected graph,
and let $S$ be a subgraph of $G$
with at least two segments such that no block of $S$ is a cycle.
Then there exists a subgraph $S'$ of $G$
 obtained from $S$ by a sequence of proper reroutings
such that every $S'$-bridge is stable.
\end{theorem}

We will need the following lemma,  a special case of
\cite[Lemma 3.2]{KawNorThoWollarge}.
A {\em linkage} in a graph is a set $\cal P$ of disjoint paths.
If  $A,B$ are sets such that each $P\in\cal P$ has one end in $A$ and the
other in $B$, then we say that $\cal P$ is a {\em linkage from $A$ to $B$}.
The {\em graph of the linkage $\cal P$} is the union of all $P\in\cal P$.
Occasionally we will use $\cal P$ in reference to the graph of $\cal P$;
thus we will use $V({\cal P})$ to denote the vertex-set of the graph of $\cal P$
and we will also speak of $\cal P$-bridges.
A path is {\em trivial} if it has exactly one vertex and {\em non-trivial}
otherwise.
By a $\cal P$-path we mean a non-trivial path with both ends in
$V({\cal P})$ and otherwise disjoint from the graph of $\cal P$.

\begin{lemma}
\mylabel{planarstrip}
Let $k\ge1$ be an integer,
let ${\cal P}=\{P_1,P_2,\ldots,P_k\}$ be a linkage in a graph $G$,
where $P_i$ has distinct ends $u_i$ and $v_i$, and let every $\cal P$-bridge
of $G$ be stable.
Assume that $G$ cannot be drawn in a disk with
$u_1,u_2,\ldots,u_k,v_k,v_{k-1},\ldots,v_1$ drawn on the boundary
of the disk in order and the paths $P_1$ and $P_k$ also drawn
on the boundary, and assume also that
there is no set $X\subseteq V(G)$ of size at most three
such that some component of $G\backslash X$ is disjoint from
$\{u_1,u_2,\ldots,u_k,v_1,v_2,\ldots,v_k\}$.
Then either
\myitem{(i)} there exist integers $i,j\in\{1,2,\ldots,k\}$ with
$|i-j|>1$ and a $\cal P$-path $Q$
in $G$ with one end in $V(P_{i})$ and
the other end in $V(P_{j})$, or
\myitem{(ii)} there exist an integer $i\in\{1,2,\ldots,k-1\}$ and two
disjoint $\cal P$-paths $Q_1$, $Q_2$ in $G$ such that $Q_j$ has ends $x_j,y_j$,
the vertices
$u_{i},x_1,x_2,v_i$ occur on $P_i$ in the order listed and
$u_{i+1},y_2,y_1,v_{i+1}$ occur on $P_{i+1}$ in the order listed, or
\myitem{(iii)} there exist an integer $i\in \{2,3,\ldots,k-1\}$ and three
$\cal P$-paths $Q_0$, $Q_{1}$, $Q_2$ such that $Q_j$ has ends $x_j$ and $y_j$,
we have
$x_0,y_0\in V(P_i)$, the vertices $x_1,x_2$ are internal vertices of $x_0P_i y_0$,
$y_1\in V(P_{i-1})$, $y_2\in V(P_{i+1})$,
and the paths $Q_0$, $Q_{1}$, $Q_2$ are pairwise disjoint, except possibly
for $x_1=x_2$.
\end{lemma}

By a {\em cylinder} we mean the surface obtained from a sphere by
removing the interiors of two disjoint closed disks $\Delta_1,\Delta_2$.
By a clockwise ordering of the boundary of $\Delta_i$ we mean
the cyclic ordering that traverses around $\Delta_i$ in clockwise
direction.
We need a slight variation of the previous lemma. We omit its proof,
because it is completely analogous.

\begin{lemma}
\mylabel{cylinderstrip}
Let $k\ge3$ be an integer,
let ${\cal P}=\{P_1,P_2,\ldots,P_k\}$ be a linkage in a graph $G$,
where $P_i$ has distinct ends $u_i$ and $v_i$, and let every $\cal P$-bridge
of $G$ be stable.
Assume that $G$ cannot be drawn in a cylinder with
$u_1,u_2,\ldots,u_k$ drawn on one boundary component in the clockwise
cyclic order listed and $v_k,v_{k-1},\ldots,v_1$ drawn on the other boundary
component  in the clockwise cyclic order listed,  assume also that
there is no set $X\subseteq V(G)$ of size at most three
such that some component of $G\backslash X$ is disjoint from
$\{u_1,u_2,\ldots,u_k,v_1,v_2,\ldots,v_k\}$,
and finally assume that if $k=3$, then no $\cal P$-bridge has vertices
of attachment on all three members of $\cal P$.
Then either
\myitem{(i)} there exist integers $i,j\in\{1,2,\ldots,k\}$ with
$|i-j|>1$ and $\{i,j\}\ne\{1,k\}$ and a $\cal P$-path $Q$
in $G$ with one end in $V(P_{i})$ and
the other end in $V(P_{j})$, or
\myitem{(ii)} there exist an integer $i\in\{1,2,\ldots,k-1\}$ and two
disjoint $\cal P$-paths $Q_1$, $Q_2$ in $G$ such that $Q_j$ has ends $x_jy_j$,
the vertices
$u_{i},x_1,x_2,v_i$ occur on $P_i$ in the order listed and
$u_{i+1},y_2,y_1,v_{i+1}$ occur on $P_{i+1}$ in the order listed, or
\myitem{(iii)} there exist an integer $i=1,2,\ldots,k$ and three
$\cal P$-paths $Q_0$, $Q_{1}$, $Q_2$ such that $Q_j$ has ends $x_j$ and $y_j$,
we have
$x_0,y_0\in V(P_i)$, the vertices $x_1,x_2$ are internal vertices of $x_0P_i y_0$,
$y_1\in V(P_{i-1})$, $y_2\in V(P_{i+1})$,
and the paths $Q_0$, $Q_{1}$, $Q_2$ are pairwise disjoint, except possibly
for $x_1=x_2$, where $P_0$ means $P_k$ and $P_{k+1}$ means $P_1$.
\end{lemma}

We finish the section by introducing several notions and a theorem from~\cite{RobSeyGM9}. We
will make use of them in the last two sections. Let $\Omega$ be a cyclic permutation of the elements of some set;
we denote this set by $V(\Omega)$. A {\em society} is a pair $(G,\Omega)$, where $G$ is a graph, and $\Omega$
is a cyclic permutation with $V(\Omega)\subseteq V(G)$. A society $(G,\Omega)$ is {\em rural} if $G$ can be drawn in a disk with $V(\Omega)$ drawn on the boundary of the disk in the order given by $\Omega$. A {\em cross} in $(G, \Omega)$ is a pair of disjoint non-trivial paths $P_1$ and $P_2$ with ends $u_1$, $v_1$ and $u_2$,$v_2$ respectively, so that $u_1, u_2, v_1, v_2 \in V(\Omega)$ appear in $\Omega$ in this or reverse order, and $P_1$ and  $P_2$ are otherwise disjoint from $V(\Omega)$.

A {\em separation} of a graph $G$ is a pair $(A,B)$ such that $A\cup B=V(G)$
and there is no edge with one end in $A-B$ and the other end in $B-A$.
The order of $(A,B)$ is $|A\cap B|$.
We say that a society $(G,\Omega)$ is {\em $4$-connected}
if there is no separation $(A,B)$ of $G$ of order at most three
with $V(\Omega)\subseteq A$ and $B-A\ne\emptyset$.

The next theorem follows from Theorems (2.3) and (2.4) in~\cite{RobSeyGM9}.
\begin{theorem}
\mylabel{thm:RSsociety}
Let $(G,\Omega)$ be a $4$-connected society with no cross. Then $(G,\Omega)$ is rural.
\end{theorem}

\section{Linear decompositions}\mylabel{sec:linear}

In this section we show that it suffices to prove  Theorem~\ref{main}
for graphs that have a ``linear decomposition" of bounded ``adhesion".
Similar techniques have been developed and used 
in~\cite{BohMahMoh,BohKawMahMoh,OpoOxlTho}.
A linear decomposition is really a tree-decomposition, where the
underlying tree is a path, but it is more convenient to number the
sets by integers rather than vertices of a path. Thus a linear
decomposition of a graph $G$ is a family of sets
${\cal W}=(W_0,W_1,\ldots,W_l)$ such that

\begin{enumerate}
\item[(L1)] $\bigcup_{i=0}^l W_i = V(G)$, and every edge of $G$ has
both ends in some $W_i$, and
\item[(L2)] if  $0\le i\le j\le k\le l$, then
$W_i \cap W_{k} \subseteq W_{j}$.
\end{enumerate}
We say that the \emph{length} of $\cal W$ is $l$.

In the proof of Theorem~\ref{main} we will need linear decompositions
satisfying the following additional properties:

\begin{enumerate}
\item[(L3)] there is an integer $q$ such that
$|W_{i-1} \cap W_{i}| = q$ for all $i = 1,2, \dots, l$,
\item[(L4)] for every $i = 1,2, \dots, l$,
$W_{i-1}\ne W_{i-1} \cap W_{i} \neq W_{i}$,
\item[(L5)] there exists a linkage
from $W_0 \cap W_1$ to $W_{l-1} \cap W_{l}$ of cardinality $q$.
\end{enumerate}

\noindent
If a linear decomposition satisfies (L3), then we say that it has
{\em adhesion} $q$.
A linkage as in (L5) will be called a {\em foundational linkage}
and its members will be called {\em foundational paths}.
We will need more properties,
but first we show that we can assume that our graph has a linear
decomposition satisfying (L1)--(L5).
In the proof we will need the following additional properties of
tree-decompositions, stated using the same notation as (W1)--(W2):

\begin{enumerate}
\item[(W3)] for every two vertices $t, t'$ of $T$ and every positive
integer $k$, either there are $k$ disjoint paths in $G$ between $Y_t$
and $Y_{t'}$, or there is a vertex $t''$
 of $T$ on the path between $t$ and $t'$ such that $|Y_{t''}| < k$,
\item[(W4)] if $t,t'$ are distinct vertices of $T$, then $Y_t \not= Y_{t'}$,
and
\item[(W5)] if $t_0 \in V(T)$ and $W$ is a component of $T - t_0$, then
$ \bigcup_{t\in V(W)} Y_t \setminus Y_{t_0} \not = \emptyset$.
\end{enumerate}

\begin{lemma}
\mylabel{lem:l5}
For all integers $k,l,p,w\ge0$ there exists an integer $N$ with
the following property.
If $G$ is a $p$-connected graph of tree-width at most $w$ with
at least $N$ vertices, then either $G$ has a minor isomorphic to
$K_{p,k}$, or $G$ has a linear decomposition of
length at least $l$ and adhesion at most $w$ satisfying {\rm (L1)--(L5)}.
\end{lemma}

\proof
Let $k,l,w\ge0$ be given integers. We will use the proof technique
of~\cite[Theorem~3.1]{OpoOxlTho} with the constants
$n_1,n_6,n_7,n_8$ and $n_9$ redefined as follows:
Let $n_1:=w$, $n_6:=l$, $n_7:=n_6^{n_1+1}$,
$n_8:={n_1\choose p}(k-1)$, and
$$n_9  := 2+n_8+n_8(n_8-1)+\cdots+n_8(n_8-1)^{\lceil n_7/2\rceil-2}.$$
We will show that $N:=n_1n_9$ satisfies the lemma.

To this end let $G$ be as stated.
The argument in Claims (1)--(4) of~\cite[Theorem~3.1]{OpoOxlTho}
shows that  $G$ either has a minor isomorphic to
$K_{p,k}$, or a tree-decomposition $(T,Y)$
satisfying (W1)--(W5) such that $T$ has a path $R$ that includes
distinct vertices $r_1,r_2,\ldots,r_l$, appearing on $R$ in the order listed,
such that for some integer $q$ with $p\le q\le w$ we have that
$|Y_{r_i}|=q$ for all $i=1,2,\ldots,l$ and $|Y_r|\ge q$ for every
$r\in V(R)$ between $r_1$ and $r_l$.

It is easy to see that there exist subtrees
$T_0,T_1,\ldots,T_l$ of $T$ such that
\myitem{(i)} $T_0\cup T_1\cup\cdots\cup T_l=T$,
\myitem{(ii)} $T_i$ and $T_j$ are disjoint whenever $|i-j|\ge2$, and
\myitem{(iii)} $V(T_{i-1})\cap V(T_i)=\{r_i\}$ for all $i=1,2,\ldots, l$.

\noindent For $i=0,1,\ldots, l$ let $W_i$ be the union of  $Y_t$
over all $t\in V(T_i)$. We claim that $(W_0,W_1,\ldots,W_l)$ is a linear decomposition of $G$
satisfying (L1)--(L5).

Property (L1) is satisfied by (W1) and (i). If  $0\le i < j < k\le l$, then for every $t \in V(T_i)$ and $t' \in V(T_k)$
the path from $t$ to $t'$ in $T$ contains the path from $r_{i+1}$ to $r_k$. Therefore,
by (W2) and (iii), we have $Y_{t} \cap Y_{t'} \subseteq Y_{r_{j}}$ and, consequently, $W_{i} \cap W_{k} \subseteq Y_{r_j} \subseteq W_j$. Thus (L2) is satisfied. Similarly, we have $W_{i-1} \cap W_i = Y_{r_i}$, and, therefore, we have $|W_{i-1} \cap W_i|=q$, implying (L3). For $1 < i \leq l$ we have $|Y_{r_{i-1}}|=|Y_{r_i}|=q$, and $Y_{r_i} \neq Y_{r_{i-1}}$ by (W4). Therefore $W_{i-1} - W_{i} \supseteq Y_{r_{i-1}} - Y_{r_i} \neq \emptyset$.  By construction,
$T_0\setminus r_1$ is the union of components of $T \setminus r_1$ disjoint from $R$. It follows from (W5) that $W_{0} - W_1 = W_0 - Y_{r_1} \neq \emptyset$.  By symmetry, $W_i - W_{i-1} \neq \emptyset$ for every $1 \leq i \leq l$, and (L4) holds. Finally, by (W3) and the choice of $r_1,r_2,\ldots,r_l$, there exists a linkage from $W_0 \cap W_1= Y_{r_1}$ to $W_{l-1} \cap W_l = Y_{r_l}$, implying (L5).
~\qed

Let $\cal P$ be a foundational linkage for a linear decomposition
${\cal W}=(W_0,W_1,\ldots,W_l)$ of a graph $G$, and let $i\in \{1,2,\ldots,l-1\}$.
We say that distinct foundational paths $P,P'\in\cal P$
are {\em bridge adjacent} in $W_i$ if there exists a $\cal P$-bridge in
$G[W_i]$ with an attachment in both $P$ and $P'$.
Given a fixed integer $p$ we
wish to consider the following properties of $\cal W$ and $\cal P$.
In our applications we will always have $p=6$.
\begin{enumerate}
\item[(L6)] for all $i\in \{1,2,\ldots,l-1\}$ and all non-trivial
paths $P\in\cal P$, if some  $\cal P$-bridge of $G[W_i]$ has
at least one attachment in $P$ and no attachment in a non-trivial
foundational path other than $P$,
then $P$ is bridge adjacent in $W_i$
to at least $p-2$ trivial members of $\cal P$,
\item[(L7)] for every $P\in \cal P$,
if there exists an index $i\in\{1,2,\ldots,l-1\}$ such that
$P[W_i]$ is a trivial path,
then $P[W_k]$ is a trivial path for all $k = 1,2, \dots, l-1$,
\item[(L8)] for every two distinct paths $P,P'\in \cal P$,
if there exists an integer $k \in \{1, \dots, l-1\}$ such that
$P$ and $P'$ are bridge adjacent in $W_k$, then they are
bridge adjacent in $W_{k'}$ for all $k' \in \{1, \dots, l-1\}$.
\end{enumerate}

With respect to condition (L8) it may be helpful to point out that
for all $i = 1,2, \dots, l$ we have $W_{i-1}\cap W_i\subseteq V(\cal P)$,
and hence each $\cal P$-bridge $H$ of $G$ satisfies $V(H)\subseteq W_k$
for some $k\in\{0,1,\ldots,l\}$, even though this index $k$ need not be unique.

\begin{lemma}
\mylabel{lem:l6}
Let $p\ge 0$ be an integer, and
let $\cal W$ be a linear decomposition of a $p$-connected graph
satisfying {\rm (L1)--(L5)}.
Then $\cal W$ has a foundational linkage $\cal P$ satisfying {\rm (L6)}.
\end{lemma}

\proof
Let ${\cal W}=(W_0,W_1,\ldots,W_l)$ be as stated.
By (L5) there exists a linkage $\cal P$ from $W_0\cap W_1$ to
$W_{l-1}\cap W_{l}$ of cardinality $q$.
Let $S$ be the union of all non-trivial paths in $\cal P$,
and let $H$ be obtained from $G[W_1\cup W_2\cup\cdots\cup W_{l-1}]$
by deleting all trivial paths in $\cal P$.
By Lemma~\ref{prestable} applied to $H$ and $S$ we may assume
(by changing $\cal P$) that  $S$
satisfies the conclusion of that lemma.
We claim that the linkage $\cal P$ then satisfies (L6).
To prove this claim suppose that
$i\in\{1,2,\ldots,l-1\}$ and some $S$-bridge $B$ of $H[W_i]$
has all its attachments in $V(P)$ for some non-trivial $P\in\cal P$;
then there are vertices $x,y\in V(P)$
such that some component $J$ of $H\backslash\{x,y\}$ has at least three vertices, includes
a vertex of $B$ and is disjoint from $V(S)-V(P)$.
Since $G$ is $p$-connected the set $V(J)$ has at least $p-2$ neighbors
among the trivial paths in $\cal P$.
Hence $P$ is bridge adjacent in $W_i$ to those trivial paths, as required.
This proves that $\cal P$ satisfies (L6).~\qed

We will make use of the following easy lemma, whose proof we omit.

\begin{lemma}
\mylabel{lem:merge}
Let ${\cal W}=(W_0,W_1,\ldots,W_l)$ be a linear decomposition
of a graph $G$ of length $l\ge2$, and let $i\in\{1,2,\ldots,l\}$.
Then
${\cal W'}:=(W_0,W_1,\ldots,W_{i-2},W_{i-1}\cup W_i,W_{i+1},W_{i+2},\ldots,W_l)$
is also a linear decomposition of $G$. Furthermore, if $\cal W$ satisfies
any of the properties {\rm (L3)--(L8)}, then so does $\cal W'$.
\end{lemma}

If $\cal W$ and $\cal W'$ are as in Lemma~\ref{lem:merge}, then we
say that $\cal W'$ was obtained from $\cal W$ by an {\em elementary
contraction}. Let $\cal P$ be a foundational linkage for $\cal W$. If $i \not \in \{1,l\}$, then let ${\cal P}':=\cal P$. If $i=1$, then let ${\cal P}'$ be the linkage obtained from $\cal P$ by restricting each $P \in \cal P$ to $W_2 \cup W_3 \cup \ldots \cup W_l$, and if $i=l$, then let ${\cal P}'$ be obtained by restricting $\cal P$ to $W_1 \cup W_2 \cup \ldots \cup W_{l-1}$. Then ${\cal P}'$ is a foundational linkage for ${\cal W}'$. It will be referred to as the \emph{corresponding restriction} of $\cal P$.
 If a linear decomposition $\cal W''$ is obtained from $\cal W$ by
a sequence of elementary contractions, then we say that $\cal W''$
is obtained from $\cal W$ by a {\em contraction}.
We will also need the following lemma about sequences of sets.

\begin{lemma}\mylabel{lem:3}
Let $l,n,\lambda\ge0 $ be  integers such that $\lambda\ge l^{n+1} n!$.
For all sequences $S_1, S_2,\dots, S_\lambda$ of subsets of $\{1, \dots, n\}$
there exist
integers $1 \le i_0 < i_1 < \dots < i_l \le \lambda+1$ such that
$$S_{i_0} \cup S_{i_0 + 1} \cup \dots \cup S_{i_1 - 1} =
S_{i_1} \cup S_{i_1 + 1} \cup \dots \cup S_{i_2 - 1} =  \dots =
S_{i_{l-1}} \cup \dots \cup S_{i_l - 1}.$$
\end{lemma}

\proof
We proceed by induction on $n$.
The lemma clearly holds when $n=0$, and so we assume that $n>0$ and
that the lemma holds for all smaller values of $n$.
If $l$ consecutive sets $S_i$ are empty, say $S_i, S_{i+1},\ldots,S_{i+l-1}$,
then the lemma holds with $i_j=i+j$ for $j=0,1,\ldots,l$.
Thus we may assume that this is not the case, and hence there
is an integer $x\in\{1,2,\ldots,n\}$ such that at least
$\lambda':=\lambda/(ln)\ge l^n (n-1)!$ of the sets $S_i$ include the element $x$.
Thus $\{1, \dots, \lambda\}$ can be partitioned into consecutive intervals
$I_1, I_2,\ldots,I_{\lambda'}$ such that each interval includes an
index $i$ with $x\in S_i$.
For $i=1,2,\ldots,\lambda'$ let $S'_i$ be the union of $S_j-\{x\}$ over
all $j\in I_i$.
By the induction hypothesis applied to the sets $S'_i$ there exist
required indices  $1\le i'_0<i'_1< \cdots <i'_l\le\lambda'+1$
for the sets $S_i'$.
For $j=0,1,\ldots,l$ let $i_j:=\min I_{i'_j}$.
It follows from the construction
that these indices satisfy the conclusion of the lemma.~\qed

\begin{lemma}
\mylabel{lem:l8}
For every triple of integers $l,p,q\ge 0$  there exists an integer
$\lambda$ such that
the following holds.
If a
graph $G$ has a linear decomposition $\cal W$ of length $\lambda+1$ and
adhesion $q$ and a foundational linkage $\cal P$
satisfying {\rm (L1)--(L6)}, then it has
a linear decomposition $\cal W'$ of length $l$ and adhesion $q$
obtained from $\cal W$ by a contraction such that $\cal W'$ and the corresponding restriction of $\cal P$
satisfy {\rm (L1)--(L8)}.
\end{lemma}

\proof
Let $l,q\ge 0$ be given, let $s:={q\choose 2}$, and let
$\mu:= l^{s+1} s!$. We will show that $\lambda:=\mu^{q+1} q!$ satisfies the
conclusion of the lemma.

Let ${\cal W}=(W_0,W_1,\ldots,W_{\lambda+1})$ be as stated.
We wish to apply Lemma~\ref{lem:3} with $q$ playing the role of $n$
and $\mu$ playing the role of $l$. For $i=1,2,\ldots,\lambda$
let $S_i$ be the set of all $P\in\cal P$
such that $P[W_i]$ is a non-trivial path.
By  Lemma~\ref{lem:3} there exist indices
$1 \le i_0 < i_1 < \dots < i_\mu \le \lambda+1$ as stated in that lemma.
Let $i_{-1}:=0$ and $i_{\mu+1}:=\lambda+1$ and
for $t=-1,0,\ldots,\mu$ define
$$W'_{t+1}:= W_{i_t}\cup W_{i_t+1}\cup\dots\cup W_{i_{t+1}-1}.$$
By Lemma~\ref{lem:merge} ${\cal W}':=(W'_0,W'_1,\ldots,W'_{\mu+1})$
is a linear decomposition of $G$ satisfying (L1)--(L6).
It follows from the construction that it also satisfies (L7).

To construct a linear decomposition satisfying (L1)--(L8) we apply
the same argument again, as follows. For a $2$-element subset
$X\subseteq\cal P$ let $S_X$ be the set of integers
$j\in\{1,2,\ldots,q\}$
such that some $\cal P$-bridge $H$ of $G$ has attachments in $P$ for
both elements $P\in X$ and satisfies $V(H)\subseteq W_j$.
By applying Lemma~\ref{lem:3} with $n:={q\choose 2}$ and $\lambda$
replaced by $\mu$
to the linear decomposition $\cal W'$
and using the same construction
we arrive at the desired linear decomposition of $G$.~\qed

Let $\cal W$ be a linear decomposition of a graph $G$
of length $l\ge2$ with foundational
linkage $\cal P$ satisfying (L1)--(L8).
We define the {\em auxiliary graph}
of the pair $({\cal W},{\cal P})$ to be the graph with vertex-set
$\cal P$ in which two paths $P,P'\in \cal P$ are adjacent if
they are bridge adjacent in $W_i$ for some (and hence every)
$i\in\{1,2,\ldots,l-1\}$.


We will need one more property of a linear decomposition $\cal W$ and its
foundational linkage ${\cal P}$.
The parameter $p$ is the same as in (L6).

\begin{enumerate}
\item[(L9)]
Let ${\cal P}_1\subseteq{\cal P}_2 \subseteq {\cal P}$ such
that $|{\cal P}_1|+|{\cal P}_2|\le p$ and  each member of ${\cal P}_1$
is non-trivial.
Then there exists a linkage $\cal Q$ in $G$ of cardinality $|{\cal P}_1|$
from $W_0\cap W_1\cap V({\cal P}_1)$ to $W_{l-1}\cap W_l\cap V({\cal P}_1)$
such that its graph is a subgraph of
$H:=G[W_0\cup W_l]\cup \bigcup_{P\in{\cal P}-{\cal P}_2} P$.
\end{enumerate}

Our objective is to show that if a graph has a linear decomposition
satisfying (L1)--(L8), then it also has one satisfying (L9).
For the proof we need a definition and a lemma.
Let ${\cal W}=(W_0,W_1,\ldots,W_l)$ be a linear decomposition
of a graph $G$ with foundational linkage ${\cal P}$ satisfying (L1)--(L8).
We say that a set ${\cal P}'$ of components of  ${\cal P}$ is {\em well-connected} if for every
two paths $P,P'\in\cal P'$ there exists a path $\cal Q$ in the auxiliary
graph of $({\cal W},{\cal P})$ such that every internal vertex of $\cal Q$
is a non-trivial foundational path belonging to ${\cal P}'$.
The lemma we need is the following.

\begin{lemma}
\mylabel{lem:prel9}
Let  $l,s,q\ge0$ be integers, and let $G$ be a graph
with a linear decomposition ${\cal W}=(W_0,W_1,\ldots,W_l)$
of length $l$, adhesion $q$
and foundational linkage ${\cal P}$ satisfying {\rm (L1)--(L8)}.
Let $\cal Q$ be a well-connected set of foundational paths,
and let $X_{ij}:= (W_{i-1}\cap W_{i}\cap V({\cal Q}))\cup
(W_j\cap W_{j+1}\cap V({\cal Q}))$.
Then for every two integers $i,j$ with $1\le i\le i+2q\le j< l$
and every two sets $A,B\subseteq X_{ij}$ of size $s$ there exist $s$
disjoint paths, each with one end in $A$,
the other end in $B$ and no internal vertex in any $W_k$
for $k\in\{0,1,\ldots,l\}-\{i,i+1,\ldots,j\}$.
\end{lemma}

\proof
Let $H$ be the subgraph of $G$ obtained by deleting $W_j-A-B$
for all  $j\in\{0,1,\ldots,l\}-\{i,i+1,\ldots,j\}$.
If the paths do not exist, then by Menger's theorem there exists
a set $Y\subseteq V(H)$ of size at most $s-1$ such that $H\backslash Y$
has no path from $A$ to $B$.
We may assume that $A\cap B=\emptyset$, for otherwise we may proceed
by induction by deleting $A\cap B$.
Since $|W_{i-1}\cap W_{i}|=|W_{j}\cap W_{j+1}|=q$ we deduce that $s\le q$.
Let $Z$ be the union of the vertex-sets of the trivial paths in $\cal P$. By (L7) and the fact that $W_i \cap W_{i+1} \subseteq V({\cal P})$ for all $i=1,2,\ldots,l-1$, the sets $W_{i+1}-Z,W_{i+3}-Z,\ldots,W_{i+2q-1}-Z$ are pairwise disjoint,
and so one of them, say $W_{m}-Z$, is disjoint from $Y$.
For $x\in X_{ij}$ let $P_x$ be the member of $\cal Q$ that includes $x$.
If $x\in W_{i-1}\cap W_i$, then let $P'_x$ denote the restriction
of $P_x$ to $W_i\cup W_{i+1}\cup\cdots\cup W_{m-1}$, and
if $x\in W_{l}\cap W_{l+1}$, then let $P'_x$ denote the restriction
of $P_x$ to $W_{m+1}\cup W_{m+2}\cup\cdots\cup W_l$.
Since the paths $P_x'$ are pairwise vertex-disjoint,
there exist $a\in A$ and $b\in B$ such that $P'_a$ and $P_b'$
are disjoint from $Y$. Since $\cal Q$ is well-connected
it follows that $P'_a\cup G[W_m]\cup P'_b$ includes a path in $H$ from
$a$ to $b$ with no internal vertex in $Z$.
That path is disjoint from $Y$, a contradiction.~\qed

\begin{lemma}\mylabel{lem:l9}
Let  $p,q\ge0$ and $l\ge3$ be integers, and let $G$ be a $p$-connected graph
with a linear decomposition
${\cal W}=(W_0,W_1,\ldots,W_{l+4q+2})$ of length $l+4q+2$, adhesion $q$
and foundational linkage $\cal P$ satisfying {\rm (L1)--(L8)}.
Let
${\cal W}':=(W'_0,W'_1,\ldots,W'_l)$,
where $W'_0:=W_0\cup W_1\cup \cdots\cup W_{2q+1}$,
$W'_i:= W_{i+2q+1}$ for $i=1,2,\ldots,l-1$ and
$W'_l:=W_{l+2q+1}\cup W_{l+2q+2}\cup \cdots\cup W_{l+4q+2}$,
and let $\cal P'$ be the corresponding restriction of $\cal P$.
Then $\cal W'$ is a linear decomposition of $G$ of length $l$ and adhesion $q$,
and $\cal P'$ is a foundational linkage for $\cal W'$
such that conditions {\rm (L1)--(L9)} hold.
\end{lemma}

\proof
The linear decomposition $\cal W'$ satisfies (L1)--(L8) by
Lemma~\ref{lem:merge}, and so it remains to show that it satisfies (L9).
Since $l\ge3$ we may choose an index $s$ with $2q+2<s<l+2q+1$.
Let ${\cal P}_1\subseteq{\cal P}_2$ be two sets of foundational paths
such that every member of ${\cal P}_1$ is non-trivial and
$|{\cal P}_1|+|{\cal P}_2|\le p$.
Let $H:=G[W_0'\cup W'_l]\cup \bigcup_{P\in {\cal P}-{\cal P}_2} P$.
We must show that there exist $|{\cal P}_1|$ disjoint paths in $H$
from $X_0:=W_0'\cap W_1'\cap V({\cal P}_1)$ to
$X_l:=W_{l-1}'\cap W_l'\cap V({\cal P}_1)$.
Since $G$ is $p$-connected and
$|W_j\cap W_{j+1}\cap V({\cal P}_2)|=|{\cal P}_2|$ we deduce
that there exists a linkage of size $|{\cal P}_1|$
from $X_0$ to $X_l$ in $G\backslash (W_s\cap W_{s+1}\cap V({\cal P}_2))$.
Let us choose such linkage, say $\cal Q$,
such that it uses the least number of edges not in $H$.
We will prove that $\cal Q$ is as desired. To do so
we may assume for a contradiction
that $\cal Q$ uses an edge
$e\in E(G)-E(H)$. By considering the linear decomposition
$(W'_l,W'_{l-1},\ldots,W'_0)$ we may assume that $e$ has both ends in
$W_i$ for some $i\in \{2q+2,2q+3,\ldots,s\}$.

By an {\em annex} we mean a maximal well-connected set
of foundational paths that includes at least one non-trivial foundational
path.
Let $\cal R$ be an annex.
We define $H_1({\cal R})$ to be the subgraph
of $J:=G[W_{1}\cup W_{2}\cup\cdots\cup W_{s}]$ consisting of the graph
of $\cal R$ restricted to $J$ and all $\cal R$-bridges that are the subgraphs of $J$ and have all vertices of attachment in $V({\cal R})$.
We define $H_0({\cal R})$ analogously as a subgraph of
$G[W_{1}\cup W_{2}\cup\cdots\cup W_{2q+1}]$.
It follows that $e$ is an edge of $H_1({\cal R})$ for some
maximal well-connected set $\cal R$ of foundational paths.
Let us assume that $e$ belongs to  $H_1({\cal R})$ for some annex $\cal R$.
Thus we fix $\cal R$ and denote $H_0({\cal R})$ and
$H_1({\cal R})$ by $H_0$ and $H_1$, respectively. We will modify the linkage $\cal Q$ within $H_1$, and will obtain
a contradiction to its choice that way.

Let $\cal Q'$ be the subset of $\cal Q$ consisting of those paths
that use at least one vertex of $H_1$. For $Q\in{\cal Q}'$ let
$a(Q)$ be its end in $X_0$, let $d(Q)$ be its end in $X_l$,
and let $b(Q)$ and $c(Q)$ be two vertices of $Q\cap H_1$ such that
the subpath of $Q$ from $b(Q)$ to $c(Q)$ is maximum and
$a(Q),b(Q),c(Q),d(Q)$ occur on $Q$ in the order listed.
It follows that $b(Q),c(Q)$ belong to $(W_0\cap W_1)\cup
(W'_0\cap W'_1)\cup (W_s\cap W_{s+1})$,
but if one of them belongs to $W'_0\cap W'_1$, then it is equal to $a(Q)$.

If $b(Q)\in W_0\cap W_1$ or $b(Q)\in W_0'\cap W_1'$ we define $b'(Q):=b(Q)$
and let $B(Q)$ be the null graph; otherwise
$b(Q)$ belongs to a foundational path $P\not \in{\cal P}_2$,
and we define $b'(Q)$ to be the unique member of
$W_{2q+1}\cap W_{2q+2}\cap V(P)$, and we let
$B(Q):=P[W_{2q+2}\cup W_{2q+3}\cup\cdots\cup W_s]$.
We define $c'(Q)$ and $C(Q)$ analogously.
By Lemma~\ref{lem:prel9} applied to $\cal W$ and $\cal P$ with
$i=0$ and $j=2q+1$
there exists a linkage $\cal S$ in $H_0$ of size $|{\cal Q}'|$ from
$\{b'(Q):Q\in{\cal Q}'\}$ to $\{c'(Q):Q\in{\cal Q}'\}$.
The fact that  $\cal R$ was chosen to be a maximal well-connected set
implies that members of this linkage are disjoint
from the members of ${\cal Q}-{\cal Q'}$.
For each $Q\in{\cal Q}'$ we delete the interior of the subpath of
$Q$ between $b(Q)$ and $c(Q)$, and add
the linkage $\cal S$ and the paths
$B(Q)$ and $C(Q)$ for all $Q\in{\cal Q}'$.
Thus we obtain a new linkage with the same properties
as $\cal Q$, but with fewer edges not in $H$, contrary to the choice
of $\cal Q$.
This completes the case when $e$ belongs to  $H_1({\cal R})$
for some annex $\cal R$, and so from now on we may assume the opposite.

Let $K$ denote the union of the trivial paths in ${\cal P}$.
Since $e$ belongs to  $H_1({\cal R})$ for no annex $\cal R$
it follows that the $K$-bridge $B$ of $H$ containing $e$ includes
no non-trivial foundational path.
Let $Q\in \cal Q$ be the path containing $e$, and let $b,c\in V(Q)$
be such that $bQc$ is a maximal subpath of $B$ containing $e$.
Since $Q$ is disjoint from $W_{s}\cap W_{s+1}\cap V({\cal P}_2)$,
and hence from the the trivial paths in ${\cal P}_2$,
we deduce that $b,c\not\in V({\cal P}_2)$.
It follows more generally
(from the fact that $e$ belongs to  $H_1({\cal R})$ for no annex $\cal R$)
that every $K$-bridge $B'$ of $H$
that has $b$ and $c$ as attachments includes no non-trivial foundational
path.
Consequently, if $B'$ includes a non-trivial subpath of some member
of $\cal Q$, then this subpath uses two vertices of $V(K)$.
On the other hand the foundational paths with vertex-sets $\{b\}$
and $\{c\}$ are adjacent in the auxiliary graph, and hence
for each $i=1,2,\ldots,q$ there exists a $K$-bridge of $G[W_i]$ whose
attachments include $b$ and $c$.
By the conclusion of the sentence before the previous one we deduce
that there is $i\in\{1,2,\ldots,q\}$ such that $W_i$ includes no
non-trivial subpath of a member of $\cal Q$.
Thus we can replace $bQc$ by a subpath of $W_i$, contrary to the choice
of $\cal Q$.
This completes the proof that $\cal W'$ and $\cal P'$ satisfy (L9).~\qed

We are now ready to state the main result of this section.

\begin{theorem}
\mylabel{thm:lindec}
For all integers $k,l,p,w\ge0$ there exists an integer $N$ with
the following property.
If $G$ is a $p$-connected graph of tree-width at most $w$ with
at least $N$ vertices, then either $G$ has a minor isomorphic to
$K_{p,k}$, or $G$ has a linear decomposition of
length at least $l$ and adhesion at most $w$ satisfying {\rm (L1)--(L9)}.
\end{theorem}

\proof Let $k,l,p,w\ge0$ be integers,
and let $l_1:=l+4w+2$.
Let $l_2$ be the minimum value of $\lambda$ such that Lemma~\ref{lem:l8}
holds for $l=l_1$, $p$ and all $q\le w$.
Finally, let $N$ be such that Lemma~\ref{lem:l5} holds for
$l=l_2$, $k,p$, and $w$. We claim that $N$ satisfies the theorem.
To prove the claim let $G$ be a $p$-connected graph of tree-width at most $w$ with
at least $N$ vertices. By Lemma~\ref{lem:l5} it has either a
minor isomorphic to $K_{p,k}$, or a linear decomposition ${\cal W}_2$
of length at least $l_2$ and adhesion $q\le w$ satisfying (L1)--(L5), and so we may assume the latter.
By Lemma~\ref{lem:l6} there is a foundational linkage ${\cal P}_1$
satisfying (L6).
By Lemma~\ref{lem:l8} the graph $G$ has a linear decomposition
${\cal W}_1$ of length $l_1$ and adhesion $q$ such that
${\cal W}_1$ and ${\cal P}_1$ satisfy (L1)--(L8).
Finally, by Lemma~\ref{lem:l9} there exist a linear decomposition
$\cal W$ of length $l$ and adhesion $q$ and a foundational linkage
satisfying (L1)--(L9).~\qed

We will need the following special case.

\begin{corollary}
\mylabel{cor:normdec}
For all integers $l,w\ge0$ there exists an integer $N$ with
the following property.
If $G$ is a $6$-connected graph of tree-width at most $w$ with
at least $N$ vertices, then either $G$ has a minor isomorphic to
$K_6$, or $G$ has a linear decomposition of
length at least $l$ and adhesion at most $w$ satisfying {\rm (L1)--(L9)} for
$p=6$.
\end{corollary}


\section{Analyzing the auxiliary graph}\mylabel{sec:auxil}

Let $G$ be a $6$-connected graph with no $K_6$ minor, and let $\cal W$
and $\cal P$ be as before and satisfy (L1)--(L9).
In this section we establish several properties of the auxiliary graph
of the pair $({\cal W}, {\cal P})$.
The first main result is Lemma~\ref{maxdeg2} stating
that if $\cal W$ is sufficiently long, then every component of
the subgraph of the auxiliary
graph induced by the non-trivial foundational paths is either a path or
a cycle.
The second main result of this section, Lemma~\ref{lem:l10},
allows us to modify the pair $({\cal W}, {\cal P})$ such that in the
new pair every non-trivial $\cal P$-bridge attaches to exactly two
non-trivial foundational paths.

\begin{figure}[h!]
\begin{center}
\includegraphics[scale = .65, page=8]{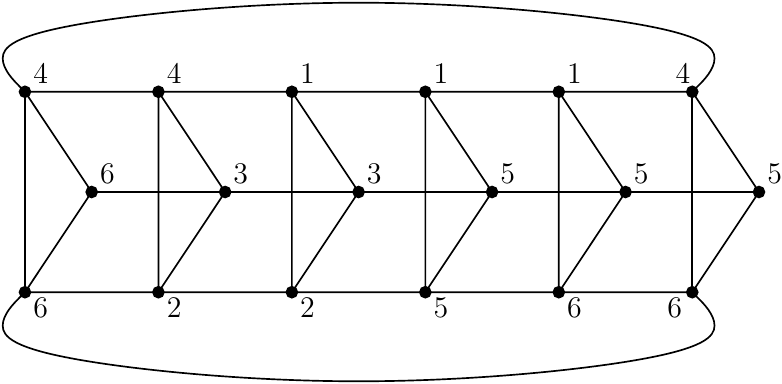}
\end{center}
\caption{Finding a $K_6$ minor in a linked $3$-cylinder of length six.}
\label{fig:3cyl}
\end{figure}

Let $k,l\ge3$ be  integers.
For $i\in \{1,2,\ldots,k\}$  let $P_i$ be a path with vertices
$v^i_1, \dots, v^i_{l}$ in order.
We define the \emph{linked k-cylinder  of length $l$} to be the graph
with vertex-set $\bigcup_{i=1}^k V(P_i)$ and edge-set
$\bigcup_{i=1}^k E(P_i) \cup \left\{  v^i_j v^{i+1}_j:  1 \le i
\le k, 1 \le j \le l \right \} \cup \{q_1, q_2\}$,
where the index notation is taken modulo $k$ and
the edges $q_1$ and $q_2$ have no common end
and each have one end in $\{ v^1_1, v^2_1, \dots, v^k_1\}$ and the other end
in $\{ v^1_{l}, v^2_{l}, \dots, v^k_l\}$.
Figure~\ref{fig:3cyl} shows a linked $3$-cylinder of length six.

\begin{lemma}\mylabel{lem:HK6}
For all integers $k \ge 3$, a linked $k$-cylinder of length twelve
has a $K_6$ minor.
\end{lemma}
\proof
By finding two suitable paths with vertex-sets in
$\{ v^i_j:  1 \le i \le k, 1 \le j \le 3 \}$, and
two paths with vertex-sets in
$\{\{ v^i_j :  1 \le i
\le k, 10 \le j \le 12 \}$, we see that a linked $k$-cylinder of length twelve
has a minor isomorphic to a linked $3$-cylinder of length six with the
additional property that the ends of the edge $q_i$
are $v_1^i$ and $v_6^i$ for $i = 1, 2$.
This graph has a $K_6$ minor as indicated in Figure~\ref{fig:3cyl}.~\qed


\begin{lemma}
\mylabel{lem:3att}
Let $l \ge 2$ and $q \ge 3$ be integers, and let ${\cal W}=(W_0,W_1,\ldots,W_l)$
be a linear decomposition of length $l$
and adhesion $q$ of a graph $G$, and let $\cal P$ be a foundational
linkage for $\cal W$ such that {\rm (L1)--(L5)} and {\rm (L9)} hold.
If for at least $48\binom{q}{3}$ indices $i\in\{1,2,\ldots,l-1\}$ there exists
a $\cal P$-bridge in $G[W_i]$ with attachments on at least three
non-trivial paths in $\cal P$, then $G$ has a $K_6$ minor.
\end{lemma}

\proof
Let $l, q$ be integers and $\mathcal{W} = (W_0, \dots, W_l)$ and $\mathcal{P}$
be given.  If there exist $48\binom{q}{3}$ distinct indices $i$ with $1 \le i \le l-1$
such that $G[W_i]$
contains a $\mathcal{P}$-bridge attaching to at least three non-trivial foundational
paths, then there exist $48$ distinct indices $i$ and three distinct non-trivial
foundational paths $P_1, P_2, P_3 \in \mathcal{P}$ such that
$G[W_i]$ contains a $\mathcal{P}$-bridge
attaching to $P_j$ for $j = 1, 2, 3$.  Then there exists a subset of
indices $I \subseteq \{1, \dots, l-1\}$ with $|I| = 24$ such that $|i - j| >2$
for all distinct $i, j \in I$, and furthermore, $G[W_i]$ contains
a bridge $B_i$ attaching to $P_j$ for all $i \in I$ and $j = 1, 2, 3$.  By property
(L9), there exist two disjoint paths $Q_1$ and $Q_2$ each with one
end in $V(P_1\cup P_2\cup P_3) \cap W_1 \cap W_2$
and one end in $V(P_1\cup P_2\cup P_3) \cap W_{l-1} \cap W_l$.
Moreover, the paths $Q_1$ and $Q_2$ do not have an internal vertex
in either $B_i \setminus V({\cal P})$ or $P_j$ for all $i \in I$ and $1 \le j \le 3$.
It follows that $G$ has a minor isomorphic to a linked $3$-cylinder of length
twelve since
each pair of successive bridges $B_i$ can be contracted to a single
cycle of length three.
By Lemma~\ref{lem:HK6} the graph $G$ has a $K_6$ minor, as desired.~\qed

The following will be a hypothesis common to several forthcoming
lemmas. In order to avoid unnecessary repetition we give it a name.

\begin{hypothesis}
\mylabel{hypot}
Let $p=6$, $l \ge 2$ and $q \ge 6$ be integers,
let $G$ be a $6$-connected graph with no $K_6$ minor,
and let ${\cal W}=(W_0,W_1,\ldots,W_l)$ be a linear decomposition
of $G$ of length $l$ and adhesion $q$ with a foundational
linkage $\cal P$ such that conditions (L1)--(L9) hold.
\end{hypothesis}

\begin{lemma}\mylabel{lem:nontriv2}
Assume Hypothesis~\ref{hypot}.  Then there do not exist
$6 \binom{q}{6}$ distinct indices $i$ with $1 \le i \le l-1$ such
that $G[W_i]$ contains a non-trivial $\mathcal{P}$-bridge attaching only
to trivial foundational paths.
\end{lemma}

\proof
Let $G$, $\cal{W}$, $\cal{P}$, $q$, and $l$ be as stated.  If
the conclusion of the lemma does not hold, then there exist
six distinct indices $i$ such that $G[W_i]$ contains a non-trivial
$\mathcal{P}$-bridge $B_i$ attaching to the same subset of
six trivial foundational paths.  By contracting the internal vertices of
each $B_i$
to a single vertex, we see $G$ would have a $K_6$ minor, a contradiction.~\qed

\begin{lemma}
\mylabel{nontrivpath}
Assume Hypothesis~\ref{hypot}. If $l> 6 \binom{q}{6}$, then $\cal P$ includes
at least one non-trivial path.
\end{lemma}

\proof
Let $G$, $\cal{W}$, $\cal{P}$, $q$, and $l$ be as stated, and suppose
for a contradiction that every path in $\mathcal{P}$ is trivial.
For every $i$, $1 \le i \le l-1$,
$G[W_i]$ contains a non-trivial bridge $B_i$, as $W_i \nsubseteq W_{i+1}$,
$W_i \nsubseteq W_{i-1}$ by (L4), in contradiction with Lemma~\ref{lem:nontriv2}.
\qed

Let $\cal W$ be a linear decomposition of a graph $G$ and let
$\cal P$ be a foundational linkage such that $\cal W$ and $\cal P$
satisfy (L1)--(L8).
By a {\em core} of the pair $({\cal W},{\cal P})$
we mean a  component of the graph obtained from
the auxiliary graph of $({\cal W},{\cal P})$ by deleting all
trivial foundational paths.
The next lemma is the first main result of this section.

\begin{lemma}
\mylabel{maxdeg2}
Assume Hypothesis~\ref{hypot}. If $l \ge 48$, then
every core of the pair $({\cal W},{\cal P})$ is a path or a cycle.
\end{lemma}

\proof Let $G$, $\cal{W}$, $\cal{P}$, $q$, and $l$ be as stated.
Suppose for a contradiction that there exists a non-trivial
foundational path $P_1 \in \cal{P}$
adjacent in the auxiliary graph to three
 non-trivial paths $P_2, P_3, P_4 \in \cal{P}$.
By property (L9), there exist two disjoint paths $Q_1$ and $Q_2$ each
with one end in $V(P_2\cup P_3\cup P_4) \cap W_0 \cap W_1$
and one end in $V(P_2\cup P_3\cup P_4)  \cap W_{l-1} \cap W_l$.
Furthermore, $Q_1$ and $Q_2$ avoid any internal vertex of $P_i$ for $1 \le i \le 4$
as well as any internal vertex of a $\mathcal{P}$-bridge in $G[W_j]$ for $1 \le j \le l-1$.
For all $i\in\{1,2,\ldots, 24\}$, we contract to a single vertex $b_i$
the set of vertices consisting
of $P_1[W_{2i-1}]$ and the
internal vertices of every non-trivial
bridge attaching to $P_1$ in $G[W_{2i-1}]$.  Note that no vertex of $Q_i$ for
$i = 1, 2$ is contained in the contracted set of $b_{2j-1}$ for any $1 \le j \le 24$.
Each vertex $b_i$ has a neighbor in each of $P_2$, $P_3$, and $P_4$.
Also, the neighbors of $b_i$ and $b_j$ are distinct for $i \neq j$.  It follows
that $G$ has a minor isomorphic to a linked
3-cylinder of length twelve, contrary to  Lemma \ref{lem:HK6}.~\qed


\begin{lemma}\mylabel{lem:4trivatt}
Assume Hypothesis~\ref{hypot}. If $l \ge 12$, then every non-trivial path in $\cal P$
is adjacent in the auxiliary graph to at most three trivial paths
in $\cal P$.
\end{lemma}

\proof  Let $G$, $\cal{W}$, $\cal{P}$, $q$, and $l$ be as stated.  Assume, to
reach a contradiction, that $P_1 \in \cal{P}$ is a non-trivial path and is
adjacent to four trivial foundational paths in the auxiliary graph.
Let the vertices comprising
the four trivial foundational paths be $v_1,v_2,v_3,v_4$.
For each $i\in\{1,2,\ldots,6\}$
we contract to a single vertex $b_i$ the vertex set containing $P_1[W_{2i-1}]$ and
the internal vertices of all non-trivial bridges of $G[W_{2i-1}]$ attaching to $P_1$.
It follows that $G$ has as a minor isomorphic to
the graph with vertex set $\{v_i: 1 \le i \le 4\}
\cup \{ b_i : 1 \le i \le 6\}$ and edges $\{v_ib_j: 1 \le i \le 4, 1 \le j \le 6\} \cup
\{b_i b_{i+1} : 1 \le i \le 5\}$.
This graph has a $K_6$ minor, and hence so does $G$, a contradiction.~\qed

\begin{corollary}\mylabel{cor:induced} Assume Hypothesis~\ref{hypot}. If $l \geq 12$, then every member of $\cal P$ is an induced path.
\end{corollary}

\proof If some non-trivial $P \in \cal P$ is not induced, then by (L6) the path $P$ is adjacent to at least $4$ trivial foundational paths in the auxilliary graph, contrary to Lemma~\ref{lem:4trivatt}.
~\qed

\begin{lemma}\mylabel{lem:3trivatt}
Assume Hypothesis \ref{hypot}.  If $l \ge 12$,  then no non-trivial foundational path
is adjacent in the auxiliary graph to three or more trivial foundational paths.
\end{lemma}

\begin{figure}[h!]
\begin{center}
\includegraphics[scale = .65, page = 12]{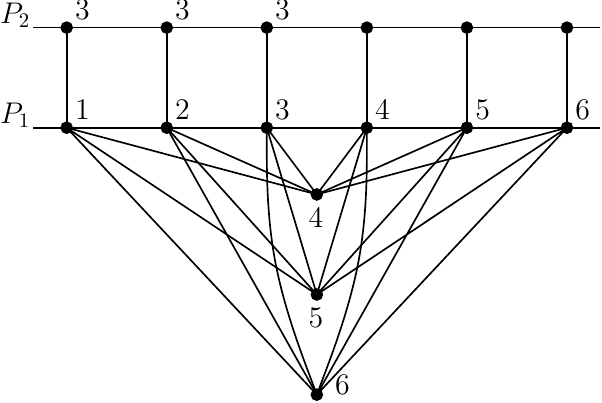}
\end{center}
\caption{Finding a $K_6$ minor when a non-trivial foundational path is bridge
adjacent to three trivial foundational paths.}
\label{fig:3trivatt}
\end{figure}

\proof Let $G$, $\cal{W}$, $\cal{P}$, $q$, and $l$ be as stated.  As above,
assume to reach a contradiction, that
$P_1 \in \cal{P}$ is a non-trivial path and is adjacent to
three trivial foundational paths in the auxiliary graph.  By the $6$-connectivity
of $G$, $P_1$ must be adjacent to another foundational path in the
auxiliary graph.  By Lemma \ref{lem:4trivatt},
such a path, call it $P_2$, must be non-trivial.
For each $i$, $1 \le i \le 6$,
we contract to a single vertex the vertex set containing $P_1[W_{2i-1}]$ and
the internal vertices of any non-trivial bridge of $G[W_{2i-1}]$ attaching to $P_1$.
It follows that $G$ has  a minor isomorphic to the graph in
Figure~\ref{fig:3trivatt}, which has a
$K_6$ minor as indicated in that figure, a contradiction.~\qed.

In the next lemma, the second main result of this section, we show that
we can assume that our linear decomposition ${\cal W}=(W_0,W_1,\ldots,W_l)$
and foundational linkage $\cal P$ satisfy the following property.
\begin{itemize}
\item[(L10)] For all $i\in\{1,2,\ldots,l-1\}$,
every non-trivial $\mathcal{P}$-bridge
of $G[W_i]$ attaches to exactly two non-trivial foundational paths.
\end{itemize}

\begin{lemma}\mylabel{lem:l10}
Assume Hypothesis~\ref{hypot}.  If
$l \ge \left ( 6\binom{q}{6} + 48 \binom{q}{3}\right )l'$,
then there exist a contraction $\mathcal{W}'$
of $\cal W$ of length $l'$ and adhesion $q$
and a foundational linkage ${\cal P}'$ for ${\cal W}'$
satisfying {\rm (L1)--(L10)}.
\end{lemma}

\proof
By Lemma \ref{lem:nontriv2} and Lemma \ref{lem:3att} and our choice of $l$,
there exists an index $\alpha$ such that for all $i\in\{1,2,\ldots,l'-1\}$,
$G[W_{\alpha + i}]$ contains neither a non-trivial $\cal P$-bridge
attaching only to trivial foundational paths nor a $\cal P$-bridge
attaching to three or more
non-trivial foundational paths.  Moreover, Lemma \ref{lem:4trivatt} and property
(L6) imply that no non-trivial bridge attaches to exactly one non-trivial
foundational path.
The lemma follows from considering the contraction
$\mathcal{W}' = \left ( \bigcup_{i = 0}^{\alpha} W_i, W_{\alpha +1},
W_{\alpha+2}, \dots, W_{\alpha + l'-1}, \bigcup_{i = \alpha + l' }^l W_i\right)$
of $\cal{W}$ and the corresponding restriction of $\cal P$.~\qed


\section{Finding and eliminating a pinwheel}\mylabel{sec:pin}

Let us assume Hypothesis~\ref{hypot}.
In the previous section we have shown that $\cal W$ and $\cal P$ can
be chosen so that for every $i\in\{1,2,\ldots,l-1\}$, every non-trivial
$\cal P$-bridge $B$ of $G[W_i]$ attaches to exactly two non-trivial
foundational paths.
The main result of this section will be used in Section~\ref{sec:taming} to show that if $G$ is not an apex graph  then
$\cal W$ and $\cal P$ can be chosen so that every such bridge attaches to no trivial foundational path.
The proof technique is different, and relies on a theory
of ``non-planar extensions" of planar graphs, developed
in~\cite{RobSeyThoExt}.

\begin{figure}[h!]
\begin{center}
\includegraphics[viewport= 0 0 400 200]{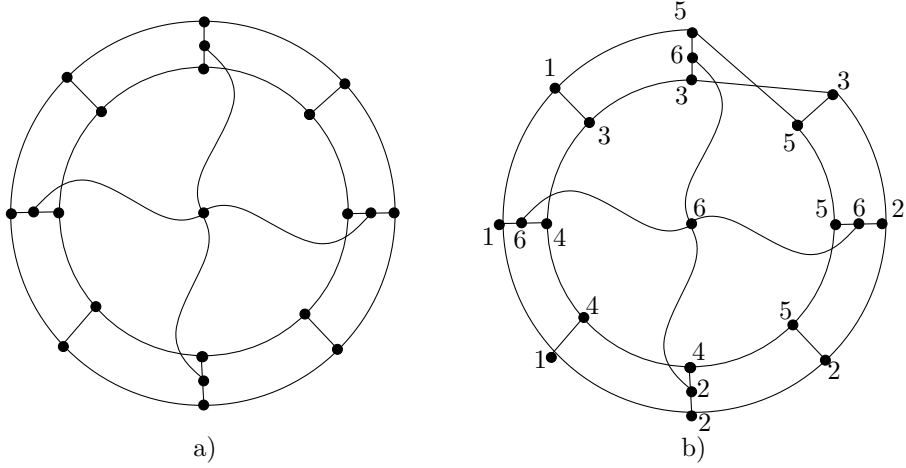}
\end{center}
\caption{(a) A pinwheel with four vanes, (b) A M\"obius pinwheel with $4$ vanes and a $K_6$ minor in it.}
\label{fig:pin}
\end{figure}

A \emph{pinwheel with $t$ vanes} is the graph defined as follows.  Let $C^1$ and $C^2$
be two disjoint cycles of length $2t$, where the vertices of $C^i$ are $v_1^i,v_2^i,\ldots,v_{2t}^i$ in order.  Let $w_1, w_2, \dots, w_t, x$
be $t+1$ distinct vertices.  The pinwheel with $t$ vanes has
vertex-set $V(C^1) \cup V(C^2) \cup \{ w_1, w_2, \dots, w_t, x\}$ and edge-set
\begin{align*}
E(C^1) &\cup E(C^2) \cup \{v^1_{2j}v^2_{2j} : 1 \le j \le t\} \\
& \cup \{ w_jv_{2j-1}^i: 1 \le j \le t, i=1,2\} \cup \{x w_j : 1 \le j \le t\}
\end{align*}
The cycles $C^1$ and $C^2$ form the \emph{rings} of the pinwheel. A pinwheel with four vanes is pictured in Figure~\ref{fig:pin}.
A \emph{M\"obius pinwheel with $t$ vanes} is obtained from a pinwheel with $t$ vanes by deleting the edges $v_{2t}^1v_1^1$ and
$v_{2t}^2v_1^2$ and adding the edges $v_{2t}^1v_1^2$ and $v_{2t}^2
v_1^1$.  The cycle formed by $V(C^1) \cup V(C^2)$ in a M\"obius pinwheel
is the \emph{ring} of the M\"obius pinwheel. A M\"obius pinwheel with $4$ vanes contains $K_6$ as a minor as shown on Figure~\ref{fig:pin}.

\begin{lemma}\mylabel{lem:pin1}
Let $q$, $l$, and $p=6$, $t\ge 4$ be positive integers.  Let $\mathcal{W} = (W_0, W_1, \dots, W_l)$
be a linear decomposition of a 6-connected graph $G$ of length $l$ and adhesion $q$ with foundational linkage ${\cal P}$
satisfying {\rm (L1)--(L9)}.  Let $P_1, P_2, P_3, Q \in {\cal P}$ be distinct, let $Q$ be trivial, and let $P_i$ be  non-trivial for $i = 1, 2, 3$.
Furthermore, let
$P_2$ be adjacent to $P_1$, $P_3$, and $Q$ in the auxiliary graph.  If
$l \ge 4t + 1$, then $G$ has  a subgraph
isomorphic to a subdivision of a pinwheel or a M\"obius pinwheel with $t$ vanes.
\end{lemma}

\proof Let $V(Q) = \{x\}$, let $P_i \cap W_0 \cap W_1 = \{s_i\}$ for $i = 1, 3$, and let $P_i \cap W_{l-1} \cap W_l = \{t_i\}$
for $i = 1, 3$.  Let $\bar {\cal P} = {\cal P} - \{P_1,P_2,P_3,Q\}$.
By property (L9),
there exist two disjoint paths $R_1$ and $R_2$ in $G[W_0 \cup W_l] \cup \bigcup_{P \in \bar{\cal P}}P$ each with one end in
$\{s_1, s_3\}$ and one end in $\{t_1, t_3\}$.  The rings of our pinwheel will be formed by
$R_1 \cup R_2 \cup P_1 \cup P_3$.  If the paths $R_1$ and $R_2$ cross, i.e. the ends of
$R_1$ are $s_1$ and $t_3$ and the ends of $R_2$ are $s_3$ and $t_1$, we construct
a M\"obius pinwheel.  Otherwise, we simply construct a pinwheel on $t$ vanes.

Note that for every $j=1,\ldots, l-1$ there exists a path $S_j$ with one end in $W_j \cap V(P_1)$ and the other end in $W_j \cap V(P_3)$, such that
$V(S_j) \subseteq W_j$, and $S_j$ is internally disjoint from
$\bigcup_{P \in {\cal P} - P_2} P$. Also, for every  $j=1,\ldots, l-1$ there exists a vertex $v_j \in W_j$ and three paths $T^1_j, T^2_j$ and $T^3_j$, internally disjoint from each other and from $\bigcup_{P \in {\cal P} - P_2} P$, satisfying the following. Each of $T^1_j, T^2_j$ and $T^3_j$ has one end $v_j$, the second end of $T^1_j$ is in $V(P_1)$, the second end of $T^3_j$ is in $V(P_3)$ and the second end of $T^2_j$ is $x$. The paths $S_j, T^1_j, T^2_j$ and $T^3_j$ are internally disjoint from the rings of our pinwheel by construction, and the paths, corresponding to the sets $W_i$ with non-consecutive indices, are also disjoint. Therefore we can use the paths corresponding to the sets $W_i$ with odd indices to construct a subgraph of $G$ isomorphic to a subdivision of a pinwheel or a M\"obius pinwheel, with  rings of the pinwheel as prescribed above.
\qed

As we have seen above a M\"obius pinwheel with sufficiently many vanes contains a $K_6$ minor. A pinwheel is, however, an apex graph.
In order to prove that graphs containing a subdivision of a pinwheel with many vanes satisfy
Theorem \ref{main}, we will need the following lemma concerning subdivisions of apex graphs
contained in larger non-apex graphs. 
The  lemma is proved in~\cite[Theorem (9.2)]{RobSeyThoExt}.

\begin{lemma}\mylabel{lem:apexsub}
Let $J$ be an internally $4$-connected triangle-free
planar graph not isomorphic to the cube, and let $F\subseteq E(J)$ be 
a nonempty set of edges such that no two edges of $F$
are incident with the same face of $J$.
Let $J'$ be obtained from $J$ by subdividing each edge in $F$
exactly once, and let
$H$ be the graph obtained from $J'$ by adding a new vertex $v\not\in V(J')$
and joining it by an edge to all the new vertices of $J'$.
Let a subdivision of $H$ be isomorphic to a subgraph of $G$, and let
$u\in V(G)$ correspond to the vertex $v$.
If $G\backslash u$ is internally $4$-connected and non-planar,
then there exists an edge $e\in E(H)$ incident with $v$ such that
either
\begin{itemize}
\item[(i)] there exist vertices $x,y\in V(J')$ not belonging
to the same face of $J'$ such that $(H\backslash e)+xy$ is
isomorphic to a minor of $G$, or
\item[(ii)] there exist vertices $x_1,x_2,x_3,x_4\in V(J')$ appearing
on some face of $J'$ in order such that $(H\backslash e)+x_1x_3+x_2x_4$ is
isomorphic to a minor of $G$.
\end{itemize}
\end{lemma}

\begin{lemma}\mylabel{lem:pin2}
If a $5$-connected graph $G$ with no $K_6$
minor contains a subdivision of a pinwheel with $20$ vanes as a subgraph, then $G$ is apex.
\end{lemma}

\proof
We will show that for every positive integer $t$ every $5$-connected non-apex graph $G$ containing a subdivision of a pinwheel with $4t$ vanes
contains a M\"obius pinwheel with $t-1$ vanes as a minor. A M\"obius pinwheel with $4$ vanes contains a $K_6$ minor, as observed above, and so the lemma will follow.

We apply Lemma~\ref{lem:apexsub}, where the graphs $H$ and $J$, the vertex $v \in V(H)$ and the set of edges $F \subseteq E(J)$ are defined as follows.
Let $H$ be the pinwheel with $4t$ vanes, and let $v$ be the ``hub"
of the pinwheel (denoted by $x$ in the definition of a pinwheel). Let the graph $J$ consist of two disjoint cycles $C^1$ and $C^2$ of length $8t$ with the vertices of $C^i = \{v_j^i: 1 \le j \le 8t\}$
for $i = 1, 2$ and $v_j^i$ adjacent to $v_{j+1}^i$ and $v_{j}^{i+1}$ for all $1 \le j \le 8t$ and $i = 1, 2$ with the
subscript addition taken modulo $8t$ and the superscript addition taken modulo $2$. Finally, let $F=\{v_{2j-1}^1v_{2j-1}^2: 1 \le j \le 4t\}$.

Suppose that outcome (ii) of Lemma~\ref{lem:apexsub} holds (the case when outcome (i) holds is analogous). If the boundary  of the face of $J$ containing the vertices
$x_1,x_2,x_3$ and $x_4$ is not one of the cycles $C_1$ and $C_2$, then without loss of generality we have $x_1=v_1^1, x_2 = v_1^2, x_3=v_2^2$ and $x_4=v_2^1$. Clearly, for every edge $e \in E(H)$ incident to $v$ the graph $(H\backslash e)+x_1x_3+x_2x_4$ contains a M\"obius pinwheel with $4t-1$ vanes as a subgraph.

Therefore, by symmetry, we assume that the vertices $x_1,x_2,x_3$ and $x_4$ are contained in $C_1$, i.e. $x_i = v_{k_i}^1$ for $i =1,2,3,4$, where, without loss of generality, $t \leq k_1, k_2,k_3, k_4 \leq 4t$. Then the subgraph $J_0$ of $J+x_1x_3+x_2x_4$ induced on $\{v_i^j : t \le i \le 4t, j=1,2\}$ contains two disjoint paths, one with ends $v_t^1$ and $v_{4t}^2$, and another with ends $v_t^2$ and $v_{4t}^1$. Now consider the graph $(H\backslash e)+x_1x_3+x_2x_4$, where $e \in E(H)$ is an edge incident to $v$, and delete all the edges of subdivision of $J_0$ from this graph, except for those that belong to the paths constructed above.
If is easy to see that the resulting graph contains a subdivision of a M\"obius pinwheel with $t-1$ vanes, as claimed.
\qed

The next corollary follows immediately from Lemmas~\ref{lem:pin1}
and~\ref{lem:pin2}.

\begin{corollary}
\mylabel{cor:mainpin}
Assume Hypothesis~\ref{hypot}.
If $l\ge81$ and some non-trivial foundational path is adjacent in
the auxiliary graph to two non-trivial and at least one trivial
foundational path, then $G$ is apex.
\end{corollary}

\section{Taming the bridges}\mylabel{sec:taming}

In Lemma~\ref{lem:l10} we have modified $\cal W$ and $\cal P$ so that
for every $i\in\{1,2,\ldots,l-1\}$
every non-trivial $\cal P$-bridge $B$ of $G[W_i]$
attaches to exactly two non-trivial
foundational paths.
Let us recall that a core is a component of the subgraph of the
auxiliary graph restricted to non-trivial foundational paths.
In this section we show that the graph consisting of all paths of a
core of $({\cal W},{\cal P})$ and all bridges that attach to two paths
of the core can be drawn in either a disk or  a cylinder, depending
on whether the core is a path or a cycle.

The following lemma follows easily from the definition of properties
(L1)--(L5) and (L9).

\begin{lemma}\mylabel{lem:reroute}
Let $l \ge 2$, $q \ge 0$, and $p \ge 0$ be integers, and let ${\cal W}=(W_0,W_1,\ldots,W_l)$
be a linear decomposition of length $l$
and adhesion $q$ of a graph $G$, and let $\cal P$ be a foundational
linkage for $\cal W$ such that {\rm (L1)--(L5)} and {\rm (L9)} hold.  
Let $i$ be fixed
with $1\le i \le l-1$ and let $Q$ be a path in $G[W_i]$ with ends $x$ and $y$
such that $x,y\in V(P)$ for some $P\in\cal P$ and
$Q$ is otherwise disjoint from $V(\cal{P})$.  
Let $P'$ be obtained from $P$ by replacing $xPy$ by $Q$.
Then the linkage $\mathcal{P}' =( \mathcal{P} -\{ P\}) \cup \{P'\}$ satisfies
{\rm (L1)--(L5)} and {\rm (L9)}.
\end{lemma}

Let $G$ be a graph and $\mathcal{W} = (W_0, \dots, W_l)$ be a
linear decomposition of length $l$ and adhesion $q$ of $G$, and let
$\cal P$ be a foundational linkage such that (L1)--(L5) hold.
Let $i\in\{1,2,\ldots,l-1\}$, let $P,P'\in\cal P$ be two non-trivial
foundational paths, let
$W_{i-1} \cap W_i \cap V(P)=\{x\}$, $W_{i-1} \cap W_i \cap V(P')=\{x'\}$,
$W_{i} \cap W_{i+1} \cap V(P)=\{y\}$, and $W_{i} \cap W_{i+1} \cap V(P')=\{y'\}$. Let
$Q_1, Q_2$ be two disjoint paths where $Q_i$ has ends $u_i$ and $v_i$ for $i=1,2$.
If the paths $Q_1$ and $Q_2$ are internally disjoint from $V(\cal{P})$,
the vertices $x$, $u_1$, $u_2$, $y$ occur
on $P$ in that order, and $x'$, $v_2$, $v_1$, $y'$ occur on $P'$ in that order,
then we say that the foundational paths $P$ and $P'$ {\em twist}.

Let $P_1$, $P_2$ and $P_3$ be three non-trivial foundational paths and let $Q_1$, $Q_2$,
and $Q_3$ be three internally disjoint paths such that $Q_j$ is also
internally disjoint from each member of $\cal{P}$ for each $j\in\{1,2,3\}$.
Let the ends of
$Q_j$ be $x_j$, $y_j$ for $1 \le j \le 3$.  The paths $Q_1$, $Q_2$, and $Q_3$
form a \emph{$P_1$-tunnel} if $x_1,y_1\in V(P_1)$,
the vertices $x_2, x_3 \in V(x_1P_1y_1) -\{x_1, y_1\}$
and $y_j \in V(P_j)$ for $j = 2, 3$.   The path $Q_1$ is called the \emph{arch}
of the tunnel.

\begin{lemma}\mylabel{lem:elimtunnels}
Let $l \ge 2$, $q \ge 3$, and $p=6$ be integers, and let ${\cal W}=(W_0,W_1,\ldots,W_l)$
be a linear decomposition of length $l$
and adhesion $q$ of a graph $G$, and let $\cal P$ be a foundational
linkage for $\cal W$ such that {\rm (L1)--(L5)} and {\rm (L9)} hold.
If there exist $48\binom{q}{3}$ distinct indices $i\in\{1,2,\ldots,l-1\}$
such that $G[W_i]$
contains a $P$-tunnel for some non-trivial foundational path $P \in \cal{P}$,
then $G$ has a $K_6$  minor.
\end{lemma}

\proof
Let $l$, $q$, $p$, $\mathcal{W}$ and $\mathcal{P}$ be given.
Assume, to reach a contradiction, that there exist $48\binom{q}{3}$ indices
$i\in\{1,2,\ldots,l-1\}$ such that $G[W_i]$ has
a  $P_i$-tunnel for some non-trivial foundational path $P_i\in\cal P$.
Reroute the  paths $P_i$ along
the arches of the $P_i$-tunnels to get a linkage $\mathcal{P}'$.
By Lemma \ref{lem:reroute} $\cal W$ and ${\cal P}'$ satisfy (L1)--(L5) and (L9).
Moreover, for each of the above $48\binom{q}{3}$ distinct indices $i$
there exists a non-trivial ${\cal P}'$-bridge in $G[W_i]$
that attaches to at least three non-trivial foundational paths.
It follows from Lemma \ref{lem:3att}
that $G$ has a  $K_6$  minor, as desired.~\qed

\begin{lemma}\mylabel{lem:elimtwists}
Let $l \ge 2$, $q \ge 3$, and $p=6$ be integers, and let $\mathcal{ W}=(W_0,W_1,\ldots,W_l)$
be a linear decomposition of length $l$
and adhesion $q$ of a graph $G$, and let $\cal P$ be a foundational
linkage for $\cal W$ such that {\rm (L1)--(L5)} and {\rm (L9)} hold.
If there exist $12\binom{q}{2}$ distinct indices $i \in \{1,2,\ldots,l-1\}$ such that $G[W_i]$ contains
a pair of twisting non-trivial foundational paths, then $G$ has a $K_6$ minor.
\end{lemma}

\begin{figure}[h!]
\begin{center}
\includegraphics[scale = .65, page = 22]{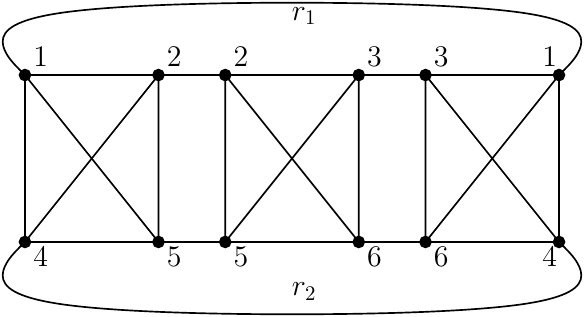}
\end{center}
\caption{Finding a $K_6$ minor when there exist a pair of non-trivial foundational
paths that twist in twelve distinct $W_i$.  The edges $r_1$ and $r_2$ are depicted
as not crossing, however, if they cross the graph still contains $K_6$ as a minor.}
\label{fig:3K4}
\end{figure}

\proof
Let $l$, $q$, $p$, $\mathcal{W}$ and $\mathcal{P}$ be given.  Assume
there exist $12\binom{q}{2}$ distinct indices $i \in \{1,2,\ldots,l-1\}$ such that $G[W_i]$ contains a pair of twisting
non-trivial foundational paths.  It follows that there exists a subset $\mathcal{I} \subseteq
\{1, 2, \dots, l-1\}$ of cardinality $12$
and non-trivial paths $P_1, P_2 \in \mathcal{P}$ such that $P_1$ and $P_2$
twist in $G[W_i]$ for all $i \in \mathcal{I}$. We use the twisting
paths to contract three disjoint $K_4$ subgraphs onto $P_1$ and $P_2$ to find
 a minor isomorphic to the graph in Figure \ref{fig:3K4}.
The edges $r_1$ and $r_2$ in the figure exist by applying property (L9)
to the ends of $P_1$ and $P_2$.
The numbering in Figure~\ref{fig:3K4} shows a $K_6$ minor,
implying that $G$ also has a $K_6$ minor, as desired.~\qed

\begin{lemma}\mylabel{lem:elimattach} Let $G$ be a $6$-connected graph with
no $K_6$ minor. Let $l \ge 2$, $q \ge 3$, and $p=6$ be integers, let $\mathcal{ W}=(W_0,W_1,\ldots,W_l)$
be a linear decomposition of length $l$
and adhesion $q$ of $G$, and let $\cal P$ be a foundational
linkage for $\cal W$ such that {\rm (L1)--(L9)} hold.
If there exist $40\binom{q}{3}$ distinct indices $i \in \{1,2,\ldots,l-1\}$ such that $G[W_i]$ contains
a non-trivial $\cal P$-bridge attaching to a trivial foundational path, then
 $G$ is apex.
\end{lemma}

\proof Let $l$, $q$, $p$, $\mathcal{W}$ and $\mathcal{P}$ be given. Assume that there exist $40\binom{q}{3}$ distinct indices $i \in \{1,2,\ldots,l-1\}$ such that $G[W_i]$ contains a non-trivial $\cal P$-bridge attaching to a trivial foundational path. By (L10) each such bridge attaches to two non-trivial foundational paths. Therefore, there exist distinct non-trivial paths $P, P' \in \cal P$ and a trivial path $Q \in \cal P$ such that $G[W_i]$ contains a $\cal P$-bridge attaching to $P,P'$ and $Q$ for at least $40$ distinct indices  $i \in \{1,2,\ldots,l-1\}$. The argument used in the proof of Lemma~\ref{lem:pin1} implies that $G$ contains a
subgraph isomorphic to a subdivision of a pinwheel with $20$ vanes or a M\"{o}bius pinwheel with $20$
vanes. Note that the M\"{o}bius pinwheel with $20$
vanes contains a $K_6$ minor, and, thus, $G$ is apex  by Lemma~\ref{lem:pin2},  as desired.~\qed

Let us assume Hypothesis~\ref{hypot}, and let $\cal C$ be a core
of $({\cal W},{\cal P})$.
We define \emph{the $i^{\hbox{th}}$ section of} $\cal C$, denoted by
$G({\cal C},i)$, to be the subgraph of $G[W_i]$, obtained from the union of
the paths in $\cal C$ and all $\cal P$-bridges of $G[W_i]$ that
attach to a member of $\cal C$ by deleting the trivial foundational paths.
By Lemma~\ref{maxdeg2} the graph $\cal C$ is a path or a cycle.
Let $P_1,P_2,\ldots,P_t$ be the vertices of $\cal C$, listed in order,
let $W_{i-1}\cap W_i\cap V(P_j)=\{u_j\}$ and
let $W_{i}\cap W_{i+1}\cap V(P_j)=\{v_j\}$.
If $\cal C$ is a path, then we say that $\cal C$ is {\em flat in $W_i$}
if $G({\cal C},i)$ can be drawn in a disk with the vertices
$u_1,u_2,\ldots,u_t,v_t,v_{t-1},\ldots,v_1$ drawn on the boundary of the
disk in order, and the paths $P_1$ and $P_t$ also drawn on the
boundary of the disk.
If $\cal C$ is a cycle, then we say that $\cal C$ is {\em flat in $W_i$}
if $G({\cal C},i)$ can be drawn in a cylinder with the vertices
$u_1,u_2,\ldots,u_t$ drawn on one of  the boundary components of the
cylinder in the clockwise order listed, and
$v_t,v_{t-1},\ldots,v_1$ drawn on the other boundary component
in the clockwise order listed.
Our next objective is to find a linear decomposition
${\cal W}=(W_0,W_1,\ldots,W_l)$ and a foundational linkage $\cal P$
such that
\begin{itemize}\item[(L11)]
Every core  of $({\cal W},{\cal P})$ is flat in $W_i$
for every $i\in\{1,2,\ldots,l-1\}$.
\item[(L12)] For every $i\in\{1,2,\ldots,l-1\}$, no  non-trivial $\cal P$-bridge of $G[W_i]$ attaches to a trivial foundational path.

\end{itemize}

\begin{lemma}\mylabel{lem:planarbridges}
Let $G$ be a 6-connected non-apex graph not containing $K_6$ as a minor.  Let $p=6$, $l \ge 2$, $q \ge 6$
be integers, and let $\mathcal{W} = (W_1, W_2, \dots, W_l)$ be a linear decomposition of $G$
of adhesion $q$ and length $l$
satisfying {\rm (L1)--(L10)}.  If $l > \left ( 88 \binom{q}{3} + 12\binom{q}{2}\right) l'$, then there exists a
contraction $\mathcal{W}'$ of $\mathcal{W}$ of length $l'$ such that ${\cal W}'$ and the corresponding restriction of $\cal P$ satisfy
{\rm (L1)--(L12)}.
\end{lemma}

\proof
Let $G$, $p$, $q$, $l$,  $\mathcal{W}$, and $\cal P$ be given.
By our choice of $l$ and Lemmas ~\ref{lem:elimtwists},~\ref{lem:elimtunnels} and~\ref{lem:elimattach},
there exists an index $\alpha$ such that
for all $i\in\{0,1,\ldots,l'\}$ the graph $G[W_{\alpha+i}]$
does not contain a $P$-tunnel for any $P$ in $\mathcal{P}$,
nor does it contain a pair of non-trivial twisting foundational paths, nor does it contain
a non-trivial bridge attaching to a trivial foundational path. We claim that the contraction
$\left ( \bigcup_{i = 0}^{\alpha-1} W_i, W_{\alpha}, W_{\alpha +1}, \dots,
W_{\alpha+l'}, \bigcup_{i = \alpha + l' +1}^l W_i\right)$ of $\mathcal{W}$ is as desired. Condition (L12) follows from the construction, and hence it suffices to prove (L11).

Fix an index $i\in\{0,1\ldots, l'\}$ and a core $\cal C$ of the auxiliary graph.
We wish to apply Lemma~\ref{planarstrip} or~\ref{cylinderstrip},
depending on whether $\cal C$ is a path or cycle, to the graph $H:=G({\cal C},\alpha + i)$
and linkage $\cal C$. Let $P_j,u_j,v_j$  for $j \in \{1,2,\ldots,t\}$ be as in the definition of flat.
By Corollary~\ref{cor:induced} and (L10) every $\cal C$-bridge of $H$ is stable, and by (L10) no $\cal C$-bridge of $H$ attaches to three or more members of $\cal C$. If there exists a set $X \subseteq V(H)$ of size at most three such that some component $J$ of $G \setminus X$ is disjoint from $\{u_1,u_2,\ldots,u_t,v_1,v_2, \ldots, v_t\}$, then by $6$-connectivity of $G$ the vertices of $J$ include a neighbor of at least three distinct trivial paths of $\cal P$. We conclude that some member of $\cal C$ is adjacent in the auxiliary graph to at least three trivial foundational paths, contrary to Lemma~\ref{lem:3trivatt}. Thus no such set $X$ exists. Next we show that none of the outcomes (i)--(iii) of Lemmas~\ref{planarstrip} and~\ref{cylinderstrip} hold. Outcome (i) does not hold by the definition of $\cal C$, and outcomes (ii) and (iii) do not hold by the choice of $\alpha$ and $i$. Thus it follows from Lemma~\ref{planarstrip} if $\cal C$ is a path or Lemma~\ref{cylinderstrip} if $\cal C$ is a cycle that $H$ can be drawn in a disk or a cylinder as described in that lemma, which is precisely the definition of $\cal C$ being flat in $W_{\alpha+i}$. Thus $\cal W'$ satisfies (L11) as well.~\qed


\section{Controlling the boundary of a planar graph}\mylabel{sec:control}

Let $G$ be a simple plane graph with the infinite region bounded by
a cycle $C$, and such that the degree of every
vertex in $V(G) - V(C)$ is at least six.
DeVos and Seymour~\cite{DevSeyExt3col} proved that
$|V(G)|\le |V(C)|^2/12 + O(|V(C)|)$.
In this section we digress to prove
 a similar result under the weaker hypothesis that $G$ has
deficiency at most five, where the {\em deficiency} of a plane graph
$G$ with the infinite region bounded by a cycle $C$
is defined
as $\sum_{v \in V(G) - V(C)} \max\{ 6- \deg(v),0\}$.
We denote the deficiency of $G$ by $\hbox{def}(G)$.
The proof is an  adaptation of the argument
from~\cite{DevSeyExt3col}, but we include it,
because the details are different.
We begin with a couple of definitions and a lemma.

A {\em quilt} is a simple plane graph $G$ with the infinite region bounded
by a cycle $C$, such that $G$ has deficiency at most five and every finite
region of $G$ is bounded by a triangle.
If exactly one vertex of $C$ has degree three, and all other vertices
have degree exactly four, then we say that $C$ is a {\em convenient graph}.
Otherwise, a {\rm convenient graph} is a subpath of $C$
with at least one edge, with both ends
of degree exactly three, and all internal vertices of degree exactly four.

\begin{lemma}\mylabel{lem:ds1}
Every quilt with no vertices of degree two has a convenient graph.
\end{lemma}

 \proof
Let $G$ be a quilt with no vertices of degree two, and
let the deficiency of $G$ be $d$. Consider the planar graph $G'$ obtained
by adding a vertex $v$ to $G$ adjacent to every vertex of $C$.
Let $|V(G)| = n$ and $|V(C)| = m$.  Then
\begin{align*} 6(n+1) - 12 & = \sum_{v \in V(G')} \deg_{G'}(v) \\
 & = \sum_{v \in V(C)} (\deg_G (v) + 1) + m + \sum_{v \in V(G) - V(C)} \deg_G(v)\\
 & \ge \sum_{v \in V(C)} \deg_G(v) + 6(n-m) - d + 2m.
\end{align*}
It follows that $\sum_{v \in V(C)} \deg_G(v) \le 4m - 6 + d$.
Since $d \le 5$ we deduce that there are strictly more vertices in $C$ of degree
three than of degree at least five.  Thus, a convenient graph exists.~\qed

The main theorem of this section follows easily from the next lemma.
If $G$ is a quilt, we define $\mu(G)$ to be $1$ if $G$ has a vertex of
degree two, and otherwise we define $\mu(G)$ to be the minimum number of
edges in a convenient graph.
Thus $\mu(G)$ is at least one, and at most the length of the cycle
bounding the infinite region of $G$.

\begin{lemma}
\mylabel{lem:quilt}
Let $G$ be a quilt on at least four vertices with the infinite region bounded
by a cycle of length $k$.
Then $|V(G)|\le k^2/2+k/2+\mu(G)+\hbox{\rm def}(G)-6$.
\end{lemma}

\proof Let $G$ and $k$ be as stated. We proceed by induction on $|V(G)|$.
If $G$ has exactly four vertices, then it is isomorphic to $K_4$, or $K_4$ minus an edge.
We have $k=3$, $\mu(G)=1$, $\hbox{def}(G)=3$, or $k=4$, $\mu(G)=1$, $\hbox{def}(G)=0$, and the lemma holds.
Thus we may assume that $G$ has at least five vertices, and that
the lemma holds for all quilts on fewer than $|V(G)|$ vertices.
 Let  $C$ be the cycle bounding the infinite region of $G$.
If $C$ has a chord, then the chord divides $G$ into two
quilts $G_1$ and $G_2$ in the obvious way.
Let the infinite region of $G_i$ have length $k_i$.
Assume first that $G_2$ has exactly three vertices.
Then by induction
\begin{align*}
|V(G)| &= |V(G_1)|+1\le k_1^2/2+k_1/2+\mu(G_1)+\hbox{def}(G_1)-6+1\\
& = k^2/2+k/2 +\mu(G_1)-k+1+\hbox{def}(G_1)-6 \\
& \le k^2/2+k/2 +\mu(G)+\hbox{def}(G)-6,
\end{align*}
as desired.
Thus we may assume that both $G_1$ and $G_2$ have at least four vertices.
Since $k_1,k_2\ge 3$ we have $3(k_1+k_2)\le k_1k_2+9$,
and hence by induction
\begin{align*}
|V(G)| &= |V(G_1)|+|V(G_2)|- 2\\
& \le  k_1^2/2+k_1/2+k_1+\hbox{def}(G_1)-6+k_2^2/2+k_2/2+k_2+\hbox{def}(G_2)-6-2\\
& = (k_1+k_2-2)^2/2  + (k_1+k_2-2)/2 +\hbox{def}(G_1)+\hbox{def}(G_2)
     -k_1k_2+3k_1+3k_2-15 \\
& \le k^2 + k/2 + \mu(G)+\hbox{def}(G)-6,
\end{align*}
as desired.
Thus we may assume that $C$ has no chord.
In particular, $G$ has no vertex of degree two.

By Lemma~\ref{lem:ds1} the quilt $G$ has a convenient graph.
Let $P$ be a convenient graph with the smallest number of edges.
Let us assume first that $P$ has exactly one edge.
Then $P$ is a path with ends $u$ and $v$, say.
Since $C$ does not have any chords and $G$ has at least five vertices,
the graph $G':=G\backslash \{u,v\}$ is a quilt.
If $G'$ has exactly three vertices, then $G$ is the wheel on five
vertices, $k=4$, $\mu(G)=1$, $\hbox{def}(G)=2$, and the lemma holds.
Thus we may assume that $G'$ has at least four vertices, and hence by
induction
\begin{align*}
|V(G)| &= |V(G')|+2\le (k-1)^2/2+(k-1)/2+\mu(G')+\hbox{def}(G')-6+2\\
& = k^2/2+k/2 +\mu(G')-k+2+\hbox{def}(G')-6 \\
& \le k^2/2+k/2 +\mu(G)+\hbox{def}(G)-6,
\end{align*}
%
%
%
as desired.
Thus we may assume that $P$ has at least two edges.
If $P=C$, then let $u$ be the unique vertex of $C$ of degree three;
otherwise $P$ is a path, and we let $u$ be an end of $P$.
Let $u'$ be the unique neighbor of $u$ that does not belong to $C$.
Then $G':=G\backslash u$ is a quilt on at least four vertices
with the infinite region bounded
by a cycle $C'$, where $C'$ has length $k$.
Since $C$ has no chords and $G$ has at least five vertices we deduce
that $\deg_{G'}(u')\ge 3$.
If equality holds, then $u$ has degree four in $G$, and hence
$\hbox{def}(G')=\hbox{def}(G)-2$.
Otherwise $\mu(G')\le\mu(G)-1$. In either case we have by induction
\begin{align*}
|V(G)| &= |V(G')|+1\le k^2/2+k/2+\mu(G')+\hbox{def}(G')-6+1\\
& \le k^2/2+k/2 +\mu(G)+\hbox{def}(G)-6,
\end{align*}
as desired.~\qed

\begin{theorem}\mylabel{thm:bdeddisc}
Let $G$ be a simple graph drawn in a disk, let $X$ be the set of vertices of $G$
drawn on the boundary of the disk, and assume that
$\sum_{v\in V(G)-X} \max\{6-\deg(v),0\}\le 5$. If $|X| \geq 3$,
then $|V(G)|\le |X|^2/2+3|X|/2-1$.
\end{theorem}

\proof Let $G$ and $X$ be as stated. We may assume, by adding edges to $G$, that $G$ is a quilt with the infinite region bounded by a cycle with vertex set $X$. By Lemma~\ref{lem:quilt} we have $|V(G)| \leq |X|^2/2 + |X|/2 + \mu(G) +\hbox{def}(G) -6 \leq |X|^2/2+3|X|/2-1$, as desired.~\qed


\section{Cylindrical tube}\mylabel{sec:cyl}

Lemma~\ref{nontrivpath} guarantees the existence of a non-empty core in a sufficiently long linear decomposition of any  $K_6$-minor-free $6$-connected
graph $G$ of bounded tree-width, assuming that such a decomposition satisfies conditions (L1)--(L9).
Lemma~\ref{maxdeg2} implies that, under the same conditions, each core is a path or a cycle. In this section we handle the case when some core of a linear decomposition of the graph $G$ is a cycle.

Before introducing the main result of this section, we need to present one more definition and a related lemma.
Let $k$, $l$ be positive integers, $k, l \ge 3$. A \emph{double crossed
$k$-cylinder of length $l$} is the graph defined as follows.  Let
$P_1, \dots, P_k$ be $k$ vertex disjoint paths with the vertex set of $P_i =
\{v^i_j: 1 \le j \le l\}$ for all $1 \le i \le k$ with $v^i_j$ adjacent to $v^i_{j+1}$
for all $1 \le j \le l-1$.  The double crossed
$k$-cylinder of length $l$ has vertex set $\{v^i_j: 1 \le j \le l, 1 \le i \le k\}$ and edge set $$\left (
 \bigcup_{i = 1}^k  E(P_i) \right) \cup \{ v^i_j v^{i+1}_j : 1 \le i \le k, 1 \le j \le l\} \cup
 \{ q_1, q_2, r_1, r_2\},$$ where the superscript addition is taken
 modulo $k$.  Furthermore, the
 ends of $q_i$ are $u_i, v_i\in \{v^j_1: 1 \le j \le  k\}$ for $i = 1, 2$ and
 the vertices $u_1, u_2, v_1, v_2$ occur in that order in the cyclic order
 $(v^1_1, v^2_1, \dots, v^1_k)$.  Similarly, the edges $r_1$ and $r_2$
 cross in the cyclic order $(v^1_l, v^2_l, \dots, v^k_l)$.  Explicitly, the ends of
 $r_i$ are $x_i, y_i \in \{v^j_l: 1 \le j \le  k\}$ for $i = 1, 2$ and occur
 in the order $x_1, x_2, y_1, y_2$ in the cyclic order
 $(v^1_l, v^2_l, \dots, v^k_l)$.

\begin{lemma}\mylabel{lem:dblextube}
Let $t$ and $l$ be integers, $t \ge 5$, $l \ge 16$.  A double crossed $t$-cylinder
of length $l$ contains $K_6$ as a minor.
\end{lemma}

\begin{figure}[h!]
\begin{center}
\includegraphics[scale = .65, page = 19]{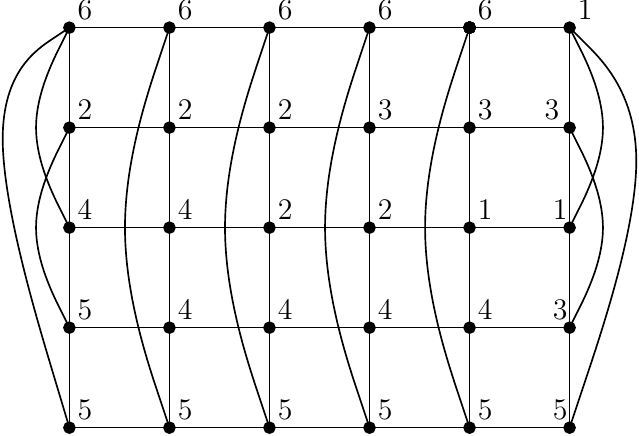}
\end{center}
\caption{A double crossed $5$-cylinder of length $6$ contains $K_6$ as a minor}
\label{fig:cyl2x}
\end{figure}

\proof
Let $G$ be a doubled crossed $t$-cylinder
of length $l$ with vertex set $\{v^i_j: 1 \le j \le l, 1 \le i \le t\}$.
By possibly routing the crossing edges $q_1$ and $q_2$ in the first five cycles
on vertices $\{v^i_j: 1 \le j \le 5, 1 \le i \le t\}$ and routing the edges $r_1$ and $r_2$ on the final
five cycles with vertex set $\{v^i_j: l-5 \le j \le l, 1 \le i \le t\}$, we see that $G$ contains as a minor
a doubled crossed $5$-cylinder $G'$
of length $6$ and moreover, with the additional property that
the ends of $q_1$ are $v^1_1$ and $v^3_1$ and the ends
of $q_2$ are $v^2_1$ and $v^4_1$.  Similarly, the edges $r_1$ and $r_2$ of $G'$
have ends $v^1_6$, $v^3_6$ and  $v^2_6$, $v^4_6$, respectively.  The graph $G$
then contains $K_6$ as a minor, as indicated in Figure \ref{fig:cyl2x}.~\qed

We now give the main result of this section.
\begin{lemma}\mylabel{lem:cylcase}
Let $p=6$, $l \ge 2$, and $q \ge 6$ be integers.
Let $G$ be a $6$-connected graph with no $K_6$ minor,
and let ${\cal W}=(W_0,W_1,\ldots,W_l)$ be a linear decomposition
of $G$ of length $l$ and adhesion $q$ with a foundational
linkage $\cal P$ satisfying {\rm(L1)--(L12)}.  Further, assume
that some core of $({\cal W},{\cal P})$ is a cycle. If
$l \ge 2q + 32$, then $G$ is apex.
\end{lemma}
\proof
Let $p$, $l$, $q$, and $\cal W$ be given, let $\cal C$ be a core of $({\cal W},{\cal P})$ that is a cycle, and assume for a contradiction that $G$ is not apex.
Let $P_1,P_2,\ldots, P_t$ be the vertices of $\cal C$ listed in order. For $i=1,2,\ldots, l-1$ let $H_i$ denote the graph $G({\cal C}, i)$, and for $j =1,2,\ldots,t$ let $u_j$ be the unique element of $V(P_j) \cap  W_q \cap W_{q+1}$ and $v_j$ the unique element of $V(P_j) \cap  W_{q+32} \cap W_{q+33}$. Let $A=\{u_1,u_2,\ldots,u_t\}$, $B=\{v_1,v_2,\ldots,v_t\}$, let $K$ denote the graph $H_{q+1} \cup H_{q+2} \cup \ldots \cup H_{q+32}$, and let $L$ denote the graph $G \setminus (V(K)-A-B)$. Since $G$ is not apex and $\cal C$ is a cycle, by Corollary~\ref{cor:mainpin} the core $\cal C$ forms a component of the auxiliary graph. Therefore, we have $K \cup L = G$ and $V(K \cap L)=A \cup B$.

We claim that $L$ does not include two disjoint paths from $A$ to $B$. Indeed, otherwise by contracting $P_i[W_{q + 2j}]$ to a single vertex for $1 \le i \le t$ and $0 \le j \le 11$,
we see that $G$ contains a linked $t$-cylinder of length twelve.
Lemma \ref{lem:HK6} then contradicts our choice of $G$. Thus there exist subgraphs $L_1,L_2$ of $L$ such that $L_1 \cup L_2=L$, $A\subseteq V(L_1)$, $B \subseteq V(L_2)$ and $|V(L_1 \cap L_2)| \leq 1$. Now property (L9) applied to $\cal C$ and a subset of $\cal C$ of size two implies that $t \geq 5$.

Let $\Omega_1$ be the cyclic permutation $(u_1,u_2,\ldots, u_t)$, and let $\Omega_2$ be the cyclic permutation $(v_1,v_2,\ldots, v_t)$. Thus $(L_1, \Omega_1)$ and
$(L_2, \Omega_2)$ are societies. Let $X= V(L_1 \cap L_2)$. By (L11) the graph $K$ can be drawn in a cylinder with $u_1,u_2,\ldots, u_t$ drawn in one boundary component in the clockwise order listed, and $v_1,v_2,\ldots, v_t$ drawn in the other boundary component in the clockwise order listed. Thus if both societies $(L_1 \setminus X, \Omega_1 \setminus X)$ and $(L_2 \setminus X, \Omega_2 \setminus X)$ are rural, then $G$ is apex, so we may assume that $(L_1 \setminus X, \Omega_1 \setminus X)$ is not rural and hence by Theorem~\ref{thm:RSsociety} it has a cross. The society  $(L_2 , \Omega_2)$ is not rural by Theorem~\ref{thm:bdeddisc}, because each vertex of $V(L_2) -  B - X$ has degree at least $6$ and $|V(L_2)| \geq qt \geq t^2= |B|^2$, because $V(L_2)$ includes each of the pairwise disjoint sets $W_i \cap W_{i+1}\cap V({\cal C})$ for $i=q+32,q+33,\ldots, 2q+31$. Likewise, $(L_2, \Omega_2)$ has a cross by Theorem~\ref{thm:bdeddisc}.

We have shown that there exist four pairwise disjoint paths, two of them forming a cross in
$(L_1, \Omega_1)$ and two forming a cross in
$(L_2, \Omega_2)$.  Let $j \in \{0,1,\ldots, 15 \}$. By the definition of core the graph $G({\cal C}, q+2j+1)$ has internally disjoint paths $Q_1,Q_2,\ldots, Q_t$ such that $Q_i$ has one end in $P_i$, the other end in $P_{i+1}$ (where $P_{t+1}$ means $P_1$),  and is otherwise disjoint from $\cal C$. Since for $j\neq j'$ the graphs $G({\cal C}, q+2j+1)$ and $G({\cal C}, q+2j'+1)$ are vertex disjoint,
we conclude that $G$ contains as a minor a double crossed $t$-cylinder
of length at least 16.  This observation contradicts Lemma \ref{lem:dblextube} and completes the
proof of the lemma.
\qed

\section{Planar strip}\mylabel{sec:planar}

We now examine the case when some core of the auxiliary graph
is a path.

\begin{lemma}\mylabel{lem:finalplan}
Let $p=6$, $l \ge 2$ and $q \ge 6$ be integers.
Let $G$ be a $6$-connected graph with no $K_6$ minor,
and let ${\cal W}=(W_0,W_1,\ldots,W_l)$ be a linear decomposition
of $G$ of length $l$ and adhesion $q$ with a foundational
linkage $\cal P$ satisfying (L1)--(L12).   Further, assume
that some core of $({\cal W},{\cal P})$ is a path. If $l \ge \max\{4q+11, 48\}$,
then $G$ is an apex graph.
\end{lemma}

\proof Let $p$, $l$, $q$, and $\cal W$ be given, let $\cal C$ be a core of $({\cal W},{\cal P})$ that is a path, and assume for a contradiction that $G$ is not apex.
Let $P_1,P_2,\ldots, P_t$ be the vertices of $\cal C$ listed in order. As in the proof of Lemma~\ref{lem:cylcase}, for $i=1,2,\ldots, l-1$ let $H_i$ denote the graph $G({\cal C}, i)$, and for $j =1,2,\ldots,t$ let $u_j$ be the unique element of $V(P_j) \cap  W_0 \cap W_{1}$ and $v_j$ the unique element of $V(P_j) \cap  W_{l-1} \cap W_{l}$. Let $A=\{u_1,u_2,\ldots,u_t\}$, $B=\{v_1,v_2,\ldots,v_t\}$, and let $\cal Q$ denote the set of trivial foundational paths adjacent in the auxiliary graph to paths in $\cal C$. Let $K$ denote the subgraph of $G$ induced on $V(H_{1} \cup H_{2} \cup \ldots \cup H_{l-1}) \cup V({\cal Q})$, and let $L$ denote the graph $G \setminus (V(K)-A-B-V({\cal Q}))$. Note that $K \cup L = G$ and $V(K) \cap V(L) = A \cup B \cup V(\cal Q)$.

We claim that either $P_1$ or $P_t$ is adjacent in the auxiliary graph to at least two paths in $\cal Q$. Suppose for a contradiction that both $P_1$ and $P_t$ are adjacent to at most one such path. We assume that $P_i$
is adjacent to exactly one trivial foundational path $S_i \in \cal Q$ for $i = 1,i=t$. The argument is similar in the case when one or both of $P_1$ and $P_t$ are not adjacent to any paths in $\cal Q$. Note that by (L12) and Corollary~\ref{cor:mainpin} all the neighbors of $V(S_1)$ and $V(S_2)$ lie on $P_1 \cup P_2$.
If $S_1 \neq S_t$, we let $\{s_i\} = V(S_i)$
for $i = 1$, $i = t$ and $K'=K$.  If $S_1 = S_t$ with $V(S_1) = V(S_t) = \{s\}$, let $K'$ be obtained from $K$ by deleting $s$, and adding new vertices $s_1$ and $s_2$, where $s_1$ is adjacent to every neighbor of $s$ on $P_1$, and
$s_t$ is adjacent to every neighbor of $s$ on $P_t$.
By property (L11), the graph
$K'$ is planar and embeds in a disk with exactly the
vertices $\{s_1, s_t\} \cup
A \cup B$ on the boundary.
Moreover, every vertex not on the boundary of the disk has degree at least six.  This is a contradiction to Theorem \ref{thm:bdeddisc},
as $|V(K')| \geq lt > (2t+2)^2$, because $l \geq 4q+11$.

\begin{figure}[h!]
\begin{center}
\includegraphics[scale = .65, page=13]{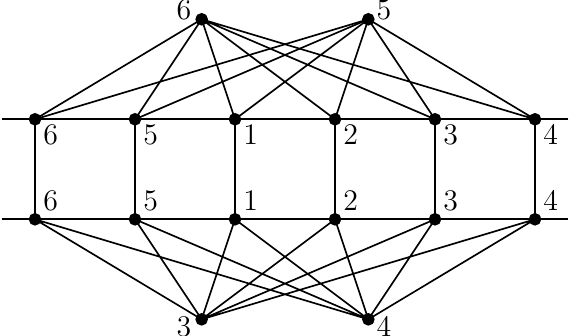}
\end{center}
\caption{Finding a $K_6$ minor when there exist four distinct trivial foundational
paths with neighbors in $\cal C$.}
\label{fig:plan2by2}
\end{figure}

\begin{figure}[h!]
\begin{center}
\includegraphics[scale = .65, page = 16]{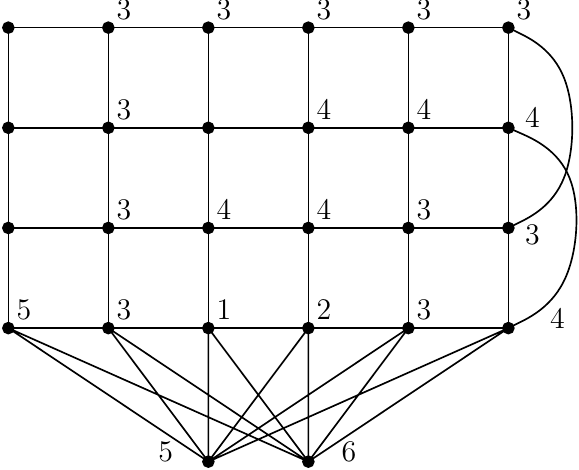}
\hfill
\includegraphics[scale = .65, page = 15]{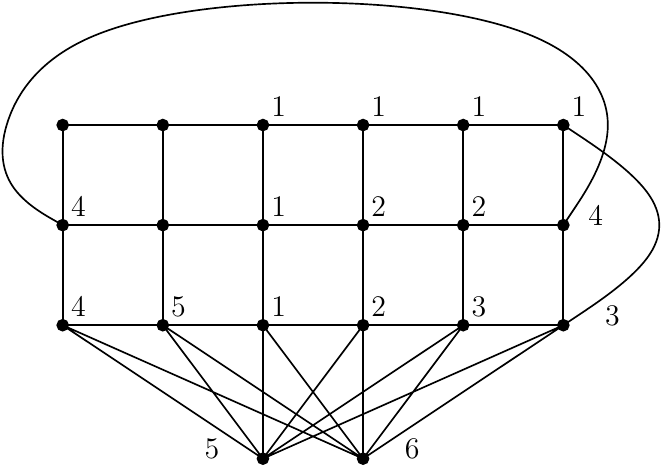}
\vskip 12pt
\vfill
\includegraphics[scale = .65, page = 14]{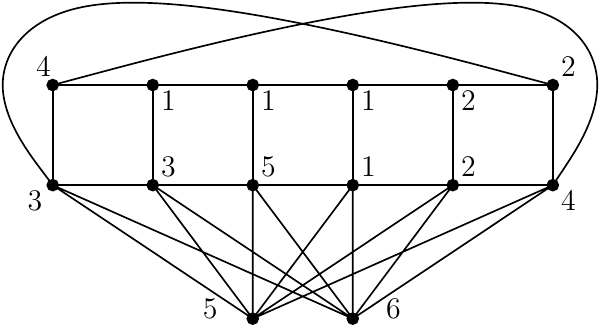}

\end{center}
\caption{Finding $K_6$ minor when the society $(L, \Omega)$ is not rural.}
\label{fig:strip3}
\end{figure}

Using the above claim and Lemma~\ref{lem:3att} we assume without loss of generality that $P_1$ is adjacent in the auxiliary graph to exactly two paths in $\cal Q$, say $Q_1$ and $Q_2$. Let $V(Q_1)=\{q_1\}$ and $V(Q_2)=\{q_2\}$. We claim that the graph $G'=G \setminus \{q_1, q_2\}$ is planar and that $P_1$ is a subset of a facial boundary of $G'$. Suppose that $P_t$ is adjacent to at least two paths in ${\cal Q} - \{Q_1,Q_2\}$.  Then  $G$ contains as a minor the graph in Figure \ref{fig:plan2by2}. The horizontal paths in the figure correspond to contractions of $P_1$ and $P_t$ and the vertical edges correspond to paths in $H_{2i+1}$ for $i=1,2,\ldots, 6$ with ends on $P_1$ and $P_t$, which exist by the definition of $\cal C$. The graph in Figure \ref{fig:plan2by2} contains a $K_6$ minor, as indicated, a contradiction.
Therefore $P_t$ is adjacent to at most one path in ${\cal Q} - \{Q_1,Q_2\}$. By  (L11), (L12) and Corollary~\ref{cor:mainpin}, the graph
$K$ is planar and embeds in the disk with $P_1$ forming part of its boundary. Let $\Omega$ be a cyclic permutation of the set $V(\Omega)= A \cup B \cup (V({\cal Q})-\{q_1,q_2\})$ ordered $u_t, u_{t-1},\ldots,u_1,v_1,\ldots,v_t$ followed by the element of $V({\cal Q})-\{q_1,q_2\}$ if $V({\cal Q})-\{q_1,q_2\} \neq \emptyset$. If the society $(L, \Omega)$ contains a cross, then $G$ contains as a minor one
of the configurations pictured in Figure \ref{fig:strip3}.  As each of this configurations contains a $K_6$ minor as indicated in Figure \ref{fig:strip3}, we conclude by Theorem~\ref{thm:RSsociety} that $(L, \Omega)$ is rural. Combined  with the planarity of $K$ this implies our claim that $G'$ is planar and $P_1$ is a subset of a facial boundary.

\begin{figure}[h!]
\begin{center}
\includegraphics[scale = .65, page = 18]{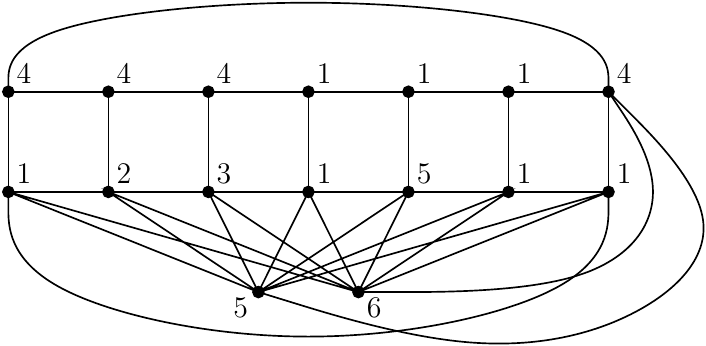}

\end{center}
\caption{Configurations giving $K_6$ minors when the trivial foundational paths
$Q_1$ and $Q_2$ have a neighbor not contained in the boundary of the face defined by $R_1 \cup P_1$}
\label{fig:plan2}
\end{figure}

Let ${\cal P}_2=\{Q_1,Q_2,P_1,P_2\}$. By property (L9), there exist two disjoint paths $R_1$ and $R_2$ in $G[W_0\cup W_l]\cup \bigcup_{P\in{\cal P}-{\cal P}_2} P$ linking the set $\{u_1, u_2\}$ to the set $\{v_1, v_2\}$.  By the claim in the previous paragraph we assume without loss of generality that $R_i$ has ends $u_i$ and $v_i$ for $i=1,2$, and that $R_1 \cup P_1$ forms a facial cycle of $G'$.  As $G$ is not apex, both $q_1$ and $q_2$ must have some neighbor not contained in $R_1 \cup P_1$. Let $q'_i$
be such a neighbor of $q_i$ for $i = 1, 2$.  The cycle $R_1 \cup P_1$ is a facial cycle in the $4$-connected planar graph $G'$, and hence there is a unique $(R_1 \cup P_1)$-bridge in $G-\{q_1,q_2\}$. It follows that for each $q'_i$ there exists a path from $q'_i$ to $R_2 \cup P_2$ avoiding $R_1 \cup P_1$.  Let $R'_i$ for $i = 1, 2$ be such paths
from $q'_i$ to $R_2 \cup P_2$.  Since $l \geq 48$ there exists an index $\alpha$ such that $W_{\alpha+i}$
is disjoint from $R'_1$ and $R'_2$ for $0 \le i \le 14$.  By considering
$P_1$ and $P_2$ and the bridges attaching to $P_1$ and $P_2$ in $H_\alpha, H_{\alpha + 1}, \ldots, H_{\alpha + 14}$, we see that
$G$ contains as a minor the graph in Figure \ref{fig:plan2}, and consequently, a $K_6$ minor, as indicated in Figure \ref{fig:plan2}. This contradiction completes the proof of the lemma.~\qed

Lemma~\ref{lem:finalplan} represents the final step in our analysis of the structure of the auxiliary graph. We are now
ready to prove Theorem~\ref{main}.

\vskip 5pt

{\noindent\bf Proof of Theorem~\ref{main}.}  Let $w\geq 1$ be an integer. Let $l_1= \max\{4w+11, 2w+32, 58\}$, let $l_2= \left (88 \binom{w}{3} + 12\binom{w}{2}\right)l_1$, and let $l_3 = \left( 6\binom{w}{6} + 48 \binom{w}{3}\right )l_2.$ By Corollary~\ref{cor:normdec} there exists an integer $N$ such that every $6$-connected graph $G$ of tree-width at most $w$ with no $K_6$ minor has a linear decomposition of length at least $l_3$ and adhesion at most $w$ satisfying properties {\rm (L1)--(L9)} for $p=6$. We claim that such an integer $N$ satisfies Theorem~\ref{main}.

Let $G$ be a $6$-connected graph of tree-width at most $w$ with at least $N$ vertices and no $K_6$ minor.
By Lemma~\ref{lem:l10} the graph $G$ has a linear decomposition of length at least $l_2$ and adhesion at most $w$ satisfying properties (L1)--(L10), and thus by Lemma~\ref{lem:planarbridges} the graph $G$ has a linear decomposition $\cal W$ of length at least $l_1$ and adhesion at most $w$ and a foundational linkage $\cal P$ satisfying properties (L1)--(L12). By Lemma~\ref{nontrivpath} $\cal P$ includes a non-trivial foundational path. By Lemma~\ref{lem:3trivatt} every non-trivial foundational path of $\cal P$ attaches to at most $2$ trivial foundational paths in the auxiliary graph. Therefore, by the $6$-connectivity of $G$, every core of $({\cal W}, {\cal P})$ has at least two vertices, and by Lemma~\ref{maxdeg2} every core is a path or a cycle. If some core of $({\cal W}, {\cal P})$ is a cycle, then $G$ is apex by Lemma~\ref{lem:cylcase}. Otherwise, $G$ is apex by Lemma~\ref{lem:finalplan}.~\qed

\section*{Acknowledgment}

We would like to acknowledge the contributions of Matthew DeVos and
Rajneesh Hegde, who worked
with us in March 2005 and contributed to this paper, 
but did not want to be included as a coauthors.

\baselineskip 11pt
\vfill
\noindent
This material is based upon work supported by the National Science Foundation.
Any opinions, findings, and conclusions or
recommendations expressed in this material are those of the authors and do
not necessarily reflect the views of the National Science Foundation.

\end{document}